\documentclass[12pt, twoside, leqno]{article}

\usepackage{amsmath,amsthm}
\usepackage{amssymb}

\usepackage{enumerate}

\usepackage{graphicx}

\usepackage{etex}

\usepackage{amsmath,amsthm,amssymb,amscd}
\usepackage{tabu}
\usepackage{amsfonts}
\usepackage{rotating}
\usepackage{multicol}

\usepackage{euscript}
\usepackage{pst-node}
\usepackage{epsfig,verbatim}
\usepackage{pictexwd,dcpic}
\usepackage[inline]{enumitem}
\usepackage{graphicx}
\usepackage{caption}
\usepackage{verbatim}

\usepackage{mathtools}

\mathtoolsset{showonlyrefs}

\def\be{\begin{equation}}
\def\ee{\end{equation}}

\def\C{{\mathbb C}} 
\def\f{\EuScript}
 
\def\P{{\mathbb P}}
\def\Z{{\mathbb Z}}
 
\def\Q{{\mathbb Q}}

\def\phi{{\varphi}}
\def\v{{\varepsilon}} 
\def\tt{\widetilde}
\def\deg{{\rm deg\,}}

\def\Gal{{\rm Gal\,}}

\def\Aut{{\rm Aut }}

\def\bp{\begin{proposition}}
\def\ep{\end{proposition}}

\def\bt{\begin{theorem}}
\def\et{\end{theorem}}
\def\br{\begin{remark}}
\def\er{\end{remark}}
\def\be{\begin{equation}}
\def\bee{\begin{equation*}}
\def\l{\label}

\def\ee{\end{equation}}
\def\eee{\end{equation*}}
\def\bl{\begin{lemma}}
\def\el{\end{lemma}}
\def\bc{\begin{corollary}}
\def\ec{\end{corollary}}
\def\pr{\noindent{\it Proof. }}

\def\bd{\begin{definition}}
\def\ed{\end{definition}}
\def\t{\widetilde}

\def\s{\underset{\mu}{\sim}}

\newtheorem{theorem}{Theorem}[section]
\newtheorem{lemma}[theorem]{Lemma}
\newtheorem{corollary}[theorem]{Corollary}
\newtheorem{proposition}[theorem]{Proposition}
\theoremstyle{definition}
\newtheorem{remark}[theorem]{Remark}

\numberwithin{equation}{section}

\begin{document}

\baselineskip=17pt


\title{On rational functions whose normalization has genus zero or one}

\author{Fedor Pakovich\\
Department of Mathematics\\ 
Ben Gurion University\\
P.O.B 653, Beer Sheva, 8410501, Israel \\
E-mail: pakovich@math.bgu.ac.il
}

\date{}

\maketitle


\renewcommand{\thefootnote}{}

\footnote{2010 \emph{Mathematics Subject Classification}: Primary 14H45, 14H30; Secondary 14G05.}

\footnote{\emph{Key words and phrases}: Galois coverings, rational points, orbifolds.}

\renewcommand{\thefootnote}{\arabic{footnote}}
\setcounter{footnote}{0}

\begin{abstract}
We give a complete list of rational functions $A$ such that the genus $g$ of the Galois closure of 
$\C(z)/\C(A)$ equals zero. We also provide a geometric description of $A$ for which $g=1.$ 

\end{abstract} 

\begin{section}{Introduction}
Let $f:\,S\rightarrow \C\P^1$ be a holomorphic function on a compact Riemann surface $S$. The normalization of $f$ is defined as a holomorphic function of the lowest possible degree
between compact Riemann surfaces $\t f:\,\t S_f\rightarrow \C\P^1$  such that $\t f$ is a Galois covering and
 $$\t f=f\circ h$$ for some  holomorphic function $h:\,\t S_f\rightarrow S$.
In this paper we study rational functions $A:\,\C\P^1\rightarrow \C\P^1$ for which the genus  of the surface $\t S_A$ equals zero or one. 
Equivalently,
we study  rational functions for which the genus of the Galois closure of the extension $\C(z)/\C(A)$ 
equals zero or one. Finally, 
the functions under consideration can be described as rational functions which are covering maps between orbifolds of non-negative Euler characteristic on the Riemann sphere.

Our main motivation for  a study of rational functions $A$ with $g(\t S_A)\leq 1$  is the fact that these functions naturally appear
in the description 
of  ``separate variable'' algebraic curves 
of genus zero and, more generally, in the  theory of functional decompositions
of rational functions.
Namely, it was shown in \cite{gen} that if an irreducible 
algebraic curve \be \l{cur} A(x)-B(y)=0,\ee where $A, B\in \C(z)$ and $\deg B\geq  \deg A,$
has genus zero, then whenever  $\deg B\geq 84\, \deg A,$
the inequality $g(\t S_A)\leq  1$ holds. 
 Moreover, for any fixed rational function $A$ with  $g(\t S_A)\leq  1$ one can find a  
rational function $B$ of arbitrary big degree such that 
corresponding curve  \eqref{cur} is irreducible and of genus zero.

Algebraic curves \eqref{cur} have been studied in  the number theory since  many interesting Diophantine equations have the form $ A(x)=B(y)$. 
By the Siegel theorem,  
an irreducible algebraic curve $\f C$ with rational coefficients may have infinitely many integer points only if $\f C$ is of genus zero with at most two points at infinity.
More generally, 
by the Faltings theorem, $\f C$ may have 
infinitely many rational points only if its genus is at most one. Therefore, 
the problem of describing curves \eqref{cur} of genus zero is important for the number theory  (see e.g. \cite{bilu}, \cite{f1},   \cite{kt}).
On the other hand,  since curve \eqref{cur} has genus zero if and only if there exist rational functions $C$ and $D$ such that 
\be \l{urik}  A\circ C=B\circ D,\ee  the problem of describing curves \eqref{cur} of genus zero is of a great value also for the decomposition theory of rational functions (see e.g. \cite{az}, \cite{prime}, \cite{alc}, \cite{amer}).

%

Another results relating  rational functions  whose normalization has genus zero or one with functional equation \eqref{urik} 
were obtained in  the paper \cite{semi} devoted to the functional  equation 
\be \l{sem} A\circ X=X\circ D,\ee   especially important for complex and arithmetic dynamics (see e.g. \cite{e2}, \cite{e}, \cite{ms}, \cite{pj}). In particular,  the results of
\cite{semi} imply that 
for  any solution $A,X,D$ of \eqref{sem} the function $X$ admits a canonical representation in the 
form 
$$X=X_0\circ W,$$ where $X_0$ satisfies $g(\t S_{X_0})\leq  1$, while $W$ is a ``compositional right factor'' of some iterate $D^{\circ k}$,
that is   $D^{\circ k}=U\circ W$ for some rational function $U$ 
(see \cite{semi}, \cite{lattes}).

In this paper we give a complete list of rational functions $A$ satisfying the condition $g(\t S_A)=0.$ Clearly, the definition implies that these functions  are exactly all possible ``compositional left factors'' of 
Galois coverings of $\C\P^1$ by $\C\P^1$. Although all such coverings were  described already by Klein, a practical calculation of their functional decompositions  is not a trivial task, and to our best knowledge 
a complete list of functions with $g(\t S_A)=0$ has never been published, although some functions from this list, and possibly even all of them, appeared here and there.     

In order to shorten the notation, 
we will 
say that rational functions $A_1$ and $A_2$ are $\mu$-equivalent and write  ${A_1}\underset{\mu}{\sim} A_2$ if 
$$A_1= \mu_1\circ A_2\circ \mu_2,$$
for some M\"obius transformations
$\mu_1$ and $\mu_2$. Our main result is  following.

\bt \l{mt1} Let $A$ be a rational function. Then $g(\t S_A)=0$ if and only if $A$ is $\mu$-equivalent to one of the functions listed below. 

\vskip 0.2cm
I) Cyclic functions:

\noindent
\begin{minipage}{1\linewidth}
\begin{equation*}
 z^n,\ \ \ n\geq 1  \leqno{a)}
\end{equation*}
\end{minipage}

\vskip 0.4cm
II) Dihedral functions: 

\noindent
\begin{minipage}{0.6\linewidth}
\begin{equation*}
 \frac{1}{2}\left(z^n+\frac{1}{z^n}\right),\ \ \ n\geq 2 \leqno{a)}
\end{equation*}
\end{minipage}%
\begin{minipage}{0.5\linewidth}
\begin{equation*}
 T_n, \ \ \ n\geq 2  \leqno{b)}
\end{equation*}
\end{minipage}

\vskip 0.5cm
III) Tetrahedral functions:

\noindent
\begin{minipage}{0.37\linewidth}
\begin{equation*}
 -\frac{1}{2^6}\frac{z^3(z^3-8)^3}{(z^3+1)^3} \leqno{a)}
\end{equation*}
\end{minipage}%
\begin{minipage}{0.3\linewidth}
\begin{equation*}
-\frac{1}{2^6}\frac{z(z-8)^3}{(z+1)^3} \leqno{b)}
\end{equation*}
\end{minipage}
\begin{minipage}{0.37\linewidth}
\begin{equation*}
 -\frac{1}{2^6}\left(\frac{z^2-4}{z-1}\right)^3  \leqno{c)}
\end{equation*}
\end{minipage}

\vskip 0.5cm
IV) Octahedral functions:

\noindent
\begin{minipage}{1\linewidth}
\begin{equation*}
 \frac{1}{2^23^3}\frac{(x^8+14x^4+1)^3}{x^4(x^4-1)^4} \leqno{a)}
\end{equation*}
\end{minipage}%

\noindent
\begin{minipage}{0.5\linewidth}
\begin{equation*}\frac{1}{2^23^3}\,{\frac { \left( {z}^{2}+14\,z+1 \right) ^{3}}{z
 \left( z-1 \right) ^{4}}}
 \leqno{b)}
\end{equation*}
\end{minipage}%
\begin{minipage}{0.5\linewidth}
\begin{equation*}
-\frac {1}{3^3}\,{\frac { \left( z^2-4 \right) ^{3} }{{z}^{
4}}}
 \leqno{c)}
\end{equation*}
\end{minipage}%

\noindent
\begin{minipage}{0.5\linewidth}
\begin{equation*}
{\frac {2^2}{3^3}}\,{\frac { \left( {z}^{4}-{z}^{2}+1 \right) ^{3}}{{z}^{
4} \left( z^2-1 \right) ^{2} }} 
 \leqno{d)}
\end{equation*}
\end{minipage}%
\begin{minipage}{0.5\linewidth}
\begin{equation*} -\frac{1}{27}\,{\frac { \left( 2\,{z}^{2}+1 \right) ^{3} \left( 2\,{z}^{2}-3
 \right) ^{3}}{ \left( 2\,{z}^{2}-1 \right) ^{4}}}
 \leqno{e)}
\end{equation*}
\end{minipage}%

\hskip -0.6cm
\begin{minipage}{0.5\linewidth}
\begin{equation*}-{\frac {2^8}{3^3}}\,{z}^{3} \left( z-1 \right) 
 \leqno{f)}
\end{equation*}
\end{minipage}%
\noindent\begin{minipage}{0.5\linewidth}
\begin{equation*}
2^8\,{\frac {z \left( z^2-7z-8 \right) ^{3}}{
 \left( {z}^{2}+20\,z-8 \right) ^{4}}}
\leqno{g)}
\end{equation*}
\end{minipage}%

\vskip 0.5cm
V) Icosahedral functions:

\noindent
\begin{minipage}{1\linewidth}
\begin{equation*} 
\frac{1}{2^63^3}\frac{
 \left( {z}^{20}+228\,{z}^{15}+494\,{z}^{10}-228\,{z}^{5}+1 \right)^{
3}} {\left( {z}^{10}-11\,{z}^{5}-1 \right) ^{5}{z}^{5}} 
\leqno{a)}
\end{equation*}
\end{minipage}%

\hskip -0.5cm
\begin{minipage}{0.6\linewidth}
\begin{equation*}-{\frac {1}{2^{11}3}}\, \left( 3\,z+5 \right) ^{3} \left( {z}^{2}+15
 \right)
 \leqno{b)}
\end{equation*}
\end{minipage}%
\begin{minipage}{0.4\linewidth}
\begin{equation*}\frac{1}{2^63^3}\frac{ \left( {z}^{2}-20 \right) ^{3}} {(z-5)}
 \leqno{c)}
\end{equation*}
\end{minipage}%

\hskip -0.5cm
\begin{minipage}{1\linewidth}
\begin{equation*}
{\frac {2^95^4}{3^2}}\,{\frac { \left( 20\,{z}^{3}-87\,z-95 \right) ^{3}
}{ \left( 20\,{z}^{2}+140\,z+101 \right) ^{5}}}\leqno{d)}
\end{equation*}
\end{minipage}%

\hskip -0.5cm
\noindent\begin{minipage}{1\linewidth}
\begin{equation*}
\frac{1}{2^63^3}\frac{
 \left( {z}^{4}+228\,{z}^{3}+494\,{z}^{2}-228\,{z}+1 \right) ^{
3}} {\left( {z}^{2}-11\,{z}-1 \right) ^{5}{z}}
 \leqno{e)}
\end{equation*}
\end{minipage}%

\hskip -0.5cm
\begin{minipage}{1\linewidth}
\begin{equation*} 
{\frac {5^4}{3^3}}\,{\frac { \left( -40\,{z}^{2}-20\,z-4 \right) ^{3}{z
}^{3} \left( 5\,{z}^{2}+5\,z+1 \right) ^{3}}{ \left( 20\,{z}^{2}+10\,z
+1 \right) ^{5}}}
\leqno{f)}
\end{equation*}
\end{minipage}%

\hskip -0.5cm
\begin{minipage}{1.1\linewidth}
\begin{equation*}
{\frac {5^3}{2^6}}\,{\frac {z \left( {z}^{2}+5\,z+40 \right) ^{3}
 \left( {z}^{2}-40\,z-5 \right) ^{3} \left( 8\,{z}^{2}-5\,z+5 \right) 
^{3}}{ \left( {z}^{4}+55\,{z}^{3}-165\,{z}^{2}-275\,z+25 \right) ^{5}}
}
\leqno{g)}
\end{equation*}
\end{minipage}%

\hskip -0.5cm
\begin{minipage}{1.1\linewidth}
\begin{equation*} 
\hskip -0.3cm  \small
\frac{1}{2^63^3}
{\frac { \left( {z}^{2}+3\,z+1 \right) ^{3}\hskip -0.1cm\left( 
{z}^{4}-4\,{z}^{3}+11\,{z}^{2}-14\,z+31 \right) ^{3}\hskip -0.1cm\left( {z}^{4}+{z
}^{3}+11\,{z}^{2}-4\,z+16 \right) ^{3}}{ \left( z-1 \right) ^{5}
 \left( {z}^{4}+{z}^{3}+6\,{z}^{2}+6\,z+11 \right) ^{5}}}
\leqno{h)}
\end{equation*} 
\end{minipage}%

\et

Notice that all the functions appeared in Theorem \ref{mt1} are Belyi functions, that is rational functions having only three critical values 
0,1, and $\infty.$ Notice also that the theorem obviously implies that any rational function $A$ of degree greater than 60  with $g(\t S_A)=0$ is either cyclic or dihedral. 
Without pretending to give a complete list of occurrences of rational functions $\mu$-equivalent to the above functions in the literature, below we point out several such examples emerging in different contexts.

The polynomials $z^n$ and $T_n$ appear in papers devoted to number theory and functional decompositions very often
(see e.g. \cite{az}, \cite{bilu},  \cite{fsch}, \cite{f1},   \cite{prime}, \cite{r2}).
In particular, the central result of the decomposition theory of polynomials, 
the so-called second Ritt theorem (see \cite{r2}), is essentially equivalent to the statement that  if  \eqref{cur} is an irreducible {\it polynomial} curve of genus zero with one point at infinity and $\deg B\geq  \deg A,$ then $A\underset{\mu}{\sim} z^n$ or $A\underset{\mu}{\sim}T_n.$
Thus, the  above mentioned result of \cite{gen} about algebraic curves \eqref{cur} can be considered as an analogue 
of the Ritt theorem for  rational functions.  

Functions ``a)'' from Theorem \ref{mt1} form a complete list of Galois coverings of $\C\P^1$ by $\C\P^1$. They were calculated in the book \cite{klein}, and nowadays
can be interpreted in terms of the ``Dessins d'enfants'' theory as Belyi functions of Platonic solids (see \cite{couv}, \cite{mazv}).
Function IV), e) is $\mu$-equivalent to
the function $3z^4-4z^3$ appearing in the paper \cite{bilu} providing a classification of  polynomial curves \eqref{cur} over $\Q$
having an infinite number of rational solutions with a bounded denominator.
Function V), b) is $\mu$-equivalent to the function $P_2$ from the paper \cite{az} about rational solutions of the functional equation 
$$A\circ C=A\circ B.$$

Functions   V), b), V), c),  and IV), e)
appear in the paper \cite{pz} about the
so-called Davenport-Zannier pairs defined over $\Q$.  
Namely, V), c) is a Belyi fun\-ction corresponding to the ``dessins $D$''  from \cite{pz} with the parameters $s=1,$ $t=1$, while V), b) and 
IV), e)  are Belyi functions corresponding to the ``dessins $A$'' with the parameters $k=s=2$, $t=1,$ and $s=3,$ $k=1,$ $t=1.$
A fun\-ction  which is $\mu$-equivalent to V), c) appears also in the paper \cite{dan} devoted to the Hall conjecture about differences between cubes and squares of integer numbers (see \cite{pz} for details).

Finally,  in the paper  \cite{mp2}, tetrahedral and octahedral functions  are used  for constructing explicit examples 
of rational functions  having decompositions into compositions of  rational functions
with a different number of indecomposable factors.

\vskip 0.2cm

In contrast to rational functions $A$ satisfying $g(\t S_A)=0$, 
functions $A$ with  $g(\t S_A)=1$ cannot be described in such an explicit way. 
Nevertheless, these functions admit quite a  precise description  in geometric terms
in the following way.

\bt \l{mt2} Let $A$ be a rational function such that $g(\t S_A)=1$. Then there exist elliptic curves 
$\f C_1$ and $\f C_2$, subgroups  $\Omega_1\subseteq {\rm Aut}({\f  C_1})$ and  $\Omega_2\subseteq {\rm Aut}({\f  C_2})$, and a holomorphic map $\alpha: \f C_1\rightarrow \f C_2$
such that the diagram 
\be \l{krotik}
\begin{CD}
\f C_1 @>\alpha>> \f C_2 \\
@VV\pi_1 V @VV\pi_2 V\\ 
\C\P^1 @>A >> \ \ \C \P^1\, ,
\end{CD}
\ee
where $\pi_1: {\f C_1}\rightarrow   
{\f C_1}/\Omega_1$ and $\pi_2: {\f C_2}\rightarrow   
{\f C_2}/\Omega_2$ are quotient maps, commutes. In the other direction, if $A$ is a rational function which makes  diagram \eqref{krotik} commutative, then $g(\t S_A)=1$, unless 
$A$ is $\mu$-equivalent either to a cyclic  function for some $n\leq 4$, or to a dihedral function for some  $n\leq 4$,  or
to a tetrahedral function. 
\et

Most known rational functions $A$ with  $g(\t S_A)=1$  are the so-called Latt\`es maps which are obtained from the above diagram for  $\f C_1=\f C_2$ and $\pi_1=\pi_2$ (see  \cite{mil2} and Section 4 below). 
Nevertheless, 
there exist  Latt\`es maps $A$ for which  $g(\t S_A)=0,$ as well as functions  $A$ with  $g(\t S_A)=1$ which are not Latt\`es maps.

The paper is organized as follows. In the second section we provide some gene\-ral definitions and results related to  functions $A$ with  $g(\t S_A)\leq 1.$ In particular, we show that such functions 
can be described in terms of their ramifications. We also give a characterization of functions $A$ with  $g(\t S_A)\leq 1$
 as covering maps between orbifolds of non-negative characteristic on the Riemann sphere. 

In the third and the fourth sections we establish some 
specific properties of rational  functions $A$ with  $g(\t S_A)=0$, and prove Theorem \ref{mt1}. We also outline a practical way of calculations of the functions from Theorem \ref{mt1} using the ``dessins d'enfants" theory.  

Finally, in the fifth section we give a  geometric characterization  of covering maps between orbifolds of zero characteristic, and 
investigate interrelations between such coverings and rational functions $A$ with  $g(\t S_A)\leq 1.$ The results of  the fifth section imply in particular Theorem \ref{mt2}. Another corollary of these results is  that for any Latt\`es map $A$ of degree greater than four the equality $g(\t S_A)=1$ holds.

\end{section}

\begin{section}{Preliminaries}
Recall that a holomorphic map $f: R_1\rightarrow R_2$  between compact  Riemann surfaces is called {\it a Galois covering} if 
its  group of deck transformations $$Aut(R_1,f)=\{h\in Aut(R_1),\ \  f\circ h=f \}$$ acts transitively on each fiber of $f.$ 
Thus, a Galois covering can be thought of as a quotient map 
\be \l{qm} R_1\rightarrow R_1/Aut(R_1,f)\cong R_2.\ee 
Equivalently, a holomorphic map  $f: R_1\rightarrow R_2$ 
is  a Galois covering if the field extension 
 $\f M(R_1)/f^{*}(\f M(R_2))$, where 
$$f^{*}:\,\f M(R_2)\rightarrow  \f M(R_1)$$ 
is the corresponding  homomorphism
of the fields of meromorphic functions, is the Galois extension. Moreover, if $f: R_1\rightarrow R_2$   a Galois covering, then 
$$\Gal(\f M(R_1)/f^{*}(\f M(R_2))\cong Aut(R_1,f)$$
(see e. g. \cite{des1}, Proposition 2.65). 
Finally, a Galois covering can be defined as a holomorphic map $f: R_1\rightarrow R_2$ such that 
\be \l{dega} \deg f=\vert Mon(f) \vert,\ee
where $Mon(f)$ is the monodromy group of $f$ (see e. g. \cite{des1}, Proposition 2.66).

Let $S$ be a compact  Riemann surface and  $f: S\rightarrow \C\P^1$  a holomorphic function. The 
{\it normalization} of $f$ is defined  as a holomorphic function of the lowest possible degree 
between compact Riemann surfaces $\t f : \t S_f\rightarrow \C\P^1$ such that $\t f$ is a Galois covering and
 $\t f=f\circ h$ for some  holomorphic function  \linebreak $h : \t S_f\rightarrow S$ (see e. g. \cite{des1}, Section 2.9).
In this paper we study rational functions $A:\C\P^1\rightarrow \C\P^1$  for which the genus  of the surface $\t S_A$ equals zero or one, or 
equivalently for which the genus of  the Galois closure of the extension $\C(z)/\C(A)$ 
equals zero or one. A convenient way for describing this class of functions in terms of their ramification 
 uses the notion of Riemann surface orbifold  (see e.g.  \cite{mil}, Appendix E, or \cite{semi}). By definition, a Riemann surface orbifold is 
a pair $\f O=(R,\nu)$ consisting of a Riemann surface $R$ 
and a ramification function $\nu:R\rightarrow \mathbb N$ which takes the value $\nu(z)=1$ except at isolated points. The Euler characteristic of an orbifold $\f O=(R,\nu)$ is defined by the formula  
\be \l{char}  \chi(\f O)=\chi(R)+\sum_{z\in R}\left(\frac{1}{\nu(z)}-1\right),\ee where $\chi(R)$ is the 
Euler characteristic of $R.$ 
For an orbifold $\f O=(R,\nu)$  we 
set $$c(\f O)=\{z_1,z_2, \dots, z_s, \dots \}=\{z\in R \mid \nu(z)>1\},$$ and 
$$\nu(\f O)=\{\nu(z_1),\nu(z_2), \dots , \nu(z_s), \dots \}.$$
For orbifolds $\f O=(R,\nu)$ and $\f O'=(R',\nu')$ we  write 
\be \l{elki} \f O\preceq \f O'\ee if $R=R'$ and for any $z\in R$ the condition $\nu(z)\mid \nu'(z)$ holds.
Clearly, \eqref{elki} implies that  $\chi(\f O) \geq \chi(\f O').$

If $R_1$, $R_2$ are Riemann surfaces provided with ramification functions $\nu_1,$ $\nu_2$, and 
$f:\, R_1\rightarrow R_2$ is a holomorphic branched covering map, then $f$
is called  {\it a covering map} $f:\,  \f O_1\rightarrow \f O_2$
{\it between orbifolds}
$\f O_1=(R_1,\nu_1)$ {\it and }$\f O_2=(R_2,\nu_2)$
if for any $z\in R_1$ the equality 
\be \l{us} \nu_{2}(f(z))=\nu_{1}(z)\deg_zf\ee holds, where $\deg_zf$ is the local degree of $f$ at the point $z$.
If $R_1$ and $R_2$ are compact and $\deg f=d,$ then 
the Riemann-Hurwitz 
formula implies that 
\be \l{rhor+} \chi(\f O_1)=d \chi(\f O_2). \ee

A universal covering of an orbifold ${\f O}$
is a covering map between orbifolds \linebreak $\theta_{\f O}:\,
\tt {\f O}\rightarrow \f O$ such that $\tt R$ is simply connected and $\tt \nu(z)\equiv 1.$ 
If $\theta_{\f O}$ is such a map, then 
there exists a group $\Gamma_{\f O}$ of conformal automorphisms of $\tt R$ such that for $z_1,z_2\in \tt R$ the equality 
$\theta_{\f O}(z_1)=\theta_{\f O}(z_2)$ holds  if and only if $z_1=\sigma(z_2)$ for some $\sigma\in \Gamma_{\f O}.$ A universal covering exists and 
is unique up to a conformal isomorphism of $\tt R,$
unless $\f O$ is the Riemann sphere with one ramified point, or  the Riemann sphere with two ramified points $z_1,$ $z_2$ such that $\nu(z_1)\neq \nu(z_2).$  Furthermore, 
$\tt R=\mathbb D$ if and only if $\chi(\f O)<0,$ $\tt R=\C$ if and only if $\chi(\f O)=0,$ and $\tt R=\C\P^1$ if and only if $\chi(\f O)>0$ (see  \cite{mil}, Appendix E, and \cite{fk}, Section IV.9.12).
Abusing  notation we will use the symbol $\tt {\f O}$ both for the
orbifold and for the  Riemann surface  $\tt R$.

For any covering map 
between orbifolds  $A:\f O_1\rightarrow \f O_2$ there exist
an isomorphism $F:\, \tt {\f O_1} \rightarrow \tt {\f O_2}$ and 
a homomorphism $\phi:\, \Gamma_{\f O_1}\rightarrow \Gamma_{\f O_2}$ such that the diagram 
\be \l{dia2}
\begin{CD}
\tt {\f O_1} @>F>> \tt {\f O_2}\\
@VV\theta_{\f O_1}V @VV\theta_{\f O_2}V\\ 
\f O_1 @>A >> \f O_2\ 
\end{CD}
\ee
commutes and 
for any $\sigma\in \Gamma_{\f O_1}$ the equality
\be \l{homm}  F\circ\sigma=\phi(\sigma)\circ F \ee holds. Vice versa, any isomorphism $F:\, \tt {\f O_1} \rightarrow \tt {\f O_2}$ satisfying \eqref{homm}  for some homomorphism $\phi:\, \Gamma_{\f O_1}\rightarrow \Gamma_{\f O_2}$ descends to a covering map $A:\f O_1\rightarrow \f O_2$  which
makes diagram \eqref{dia2} commutative (see e.g, \cite{semi}, Proposition 3.1).

With each holomorphic function $f:\, R_1\rightarrow R_2$ between compact Riemann surfaces 
one can associate in a natural way two orbifolds $\f O_1^f=(R_1,\nu_1^f)$ and 
$\f O_2^f=(R_2,\nu_2^f)$, setting $\nu_2^f(z)$  
equal to the least common multiple of the local degrees of $f$ at the points 
of the preimage $f^{-1}\{z\}$, and $$\nu_1^f(z)=\nu_2^f(f(z))/\deg_zf.$$ We will call $\f O_2^f$ {\it the ramification orbifold} of $f$.
By construction,  \be \l{covve} f:\, \f O_1^f\rightarrow \f O_2^f\ee
is a covering map between orbifolds.  
Furthermore, 
it is easy to  see that the covering map $f:\, \f O_1^f\rightarrow \f O_2^f$ is minimal in the following sense. For any covering map between orbifolds $f:\, \f O_1\rightarrow \f O_2$ we have:
\be \l{elki+} \f O_1^f\preceq \f O_1, \ \ \ \f O_2^f\preceq \f O_2.\ee
Orbifolds $\f O_1^f$ and $\f O_2^f$ always have a universal covering
(see \cite{semi}, Lemma 4.2).

Since any Galois covering $f: R_1\rightarrow R_2$ is quotient map \eqref{qm},
for any branch point $z_i$, $1\leq i \leq r,$  of $f$ there exists a number $d_i$ such that $f^{-1}\{z_i\}$ consists of $\vert Aut(R_1,f)\vert /d_i$ points, and at each of these points the multiplicity of $f$ equals $d_i$.
Indeed, points of $f^{-1}\{z_i\}$ form a single orbit of $Aut(R_1,f)$ and have conjugated stabilizers of the same order. 
In the above notation, we can formulate this property of Galois coverings as follows: for any Galois covering $f: R_1\rightarrow R_2$
the orbifold $\f O_1^f$ is non-ramified.

The following statement coincides with Lemma 1 of  \cite{gen}.
For the reader's convenience we repeat the arguments.

\bl \l{ml} 
Let $A$ be a  rational function. Then $g(\t S_A)=0$ if and only if  $\chi(\f O_2^A)> 0$, and  $g(\t S_A)=1$ if and only if  $\chi(\f O_2^A)= 0$. 
\el
\pr 
Let $f: S\rightarrow \C\P^1$ be a Galois covering of $\C\P^1$. 
Applying the Riemann-Hurwitz formula, we see that 
$$2g(S)-2=-2\vert \Gamma\vert +\sum_{i=1}^r\frac{\vert \Gamma\vert}{d_i}\left(d_i-1\right),$$ where $\Gamma=Aut(S,f)$, 
implying that
\be \l{gc} \chi(\f O_2^f)=2+\sum_{i=1}^r\left(\frac{1}{d_i}-1\right)=\frac{2-2g(S)}{\vert \Gamma\vert }.\ee 
Thus, if $f: S\rightarrow \C\P^1$ is a Galois covering, then $g(S)=0$ if and only if $\chi(\f O_2^f)>0$, while  $g(S)=1$ if and only if $\chi(\f O_2^f)=0$.

Let now $A:\, \C\P^1\rightarrow \C\P^1$ be an arbitrary rational function of degree $n$. Since the normalization
$\t A:\, \t S_A\rightarrow \C\P^1$ of $A$ can be described as any  connected
component of the $n$-fold  fiber product of $A$ distinct from the diagonal components (see \cite{fried}, $\S$I.G), it follows from 
the construction of the fiber product (see e. g. \cite{prime}, Section 2 and 3) that 
\be \l{vtgh} \f O_2^A=\f O_2^{\t A}.\ee  Thus, $g(\t S_A)=0$ if and only if  $\chi(\f O_2^A)> 0$, and  $g(\t S_A)=1$ if and only if  $\chi(\f O_2^A)= 0$. 
\qed

\vskip 0.2cm

Lemma \ref{ml} gives a simple practical way for checking whether $g(\t S_A)\leq 1$ in terms of the ramification of $A$.

\bc \l{cor1}  Let $A$ be a rational function. Then  $g(\t S_A)=0$ if and only if   $\nu(\f O_2^A)$ belongs to the  list  
\be \l{list} \{2,3,6\}, \ \ \ \{2,4,4\}, \ \ \ \{3,3,3\}, \ \ \ \{2,2,2,2\},\ee
while  $g(\t S_A)>0$ if and only if  $\nu(\f O_2^A)$ belongs to the  list 
 \be \l{list2} \{n,n\}, \ \ n\geq 1,  \ \ \ \{2,2,n\}, \ \ n\geq 2,  \ \ \ \{2,3,3\}, \ \ \ \{2,3,4\}, \ \ \ \{2,3,5\}.\ee
\ec
\pr Indeed, it is well known and follows easily from \eqref{char} that if  $\f O$ is an orbifold on $\C\P^1$ having a universal covering, then  $\chi(\f O)=0$  if and only if $\nu(\f O)$  belongs to list  \eqref{list}, and  $\chi(\f O)>0$ if and only if $\nu(\f O)$  belongs to list \eqref{list2}. \qed
\vskip 0.2cm

Another corollary of Lemma \ref{ml} is the following statement. 

\bc \l{cor2} Let $A$ be a rational function. Then  $g(\t S_A)=0$ if and only if  there exist 
orbifolds of positive Euler characteristic $\f O_1$ and $\f O_2$ on $\C\P^1$ such that  $A:\, \f O_1\rightarrow \f O_2$ is a covering map between orbifolds. 
Similarly, if $g(\t S_A)=1$, then   there exist 
orbifolds of zero Euler characteristic such that $A:\, \f O_1\rightarrow \f O_2$ is a covering map. On the other hand, the fact that 
$A:\, \f O_1\rightarrow \f O_2$ is a covering map between orbifolds of zero  Euler characteristic implies only that 
$g(\t S_A)\leq 1$.
\ec
\pr 
Indeed, if $\chi(\f O_2^A)> 0,$ then \eqref{rhor+} implies that  $\chi(\f O_1^A)> 0$ and hence  $A:\f O_1^A\rightarrow \f O_2^A$  is a covering map between  orbifolds of positive Euler characteristic.
Similarly, $\chi(\f O_2^A)= 0$ implies that  $\chi(\f O_1^A)= 0$.

In the other direction, if $A: \f  O_1\rightarrow \f O_2$ is a covering map between some orbifolds of positive Euler characteristic, then \eqref{elki+} implies that  
$\f O_1^A$ and $\f O_2^A$ also have positive Euler characteristic. On the other hand, if $\f O_1$ and $\f O_2$ have zero Euler characteristic, then  \eqref{elki+} implies only that  
$\f O_1^A$ and $\f O_2^A$ have non-negative Euler characteristic. \qed

\end{section}

\begin{section}{Functions with $\chi(\f O_2^A)>0$}

Let $f,g$ be rational functions. We will call $g$ a {\it compositional left factor} of $f$ if 
$f=g\circ h$ for some rational function $h.$
It is clear that rational functions $A$ with $g(\t S_A)= 0$
are exactly compositional left factors of Galois coverings  $f:\, \C\P^1\rightarrow \C\P^1.$ 
Notice that since for a Galois covering  $f:\, \C\P^1\rightarrow \C\P^1$ 
the orbifold  $\f O_1^f$ is non-ramified, any Galois covering  $f:\, \C\P^1\rightarrow \C\P^1$  is a universal covering of the orbifold   $\f O=\f O_2^f$ with $\chi(\f O)>0$, and vice versa 
for any orbifold $\f O$ with $\chi(\f O)>0$  its universal covering 
$$\theta_{\f O}:\, \C\P^1\rightarrow \C\P^1/\Gamma_{\f O} \cong \C\P^1$$ is 
a Galois covering. Thus, 
we can identify Galois coverings  $f:\, \C\P^1\rightarrow \C\P^1$ with universal coverings of orbifolds $\f O$ of positive Euler characteristic on $\C\P^1$.

Recall that any finite subgroup of $ Aut(\C\P^1) $ is
isomorphic to one of the following groups:
\be \l{wsx} \Z/n\Z,\  n\geq 1, \ \ \ D_{2n}, \  n\geq 2,  \ \ \  A_4, \ \ \ S_4, \ \ \ A_5.\ee 
Moreover,  these isomorphism classes are also conjugacy classes. 
Groups \eqref{wsx} are groups $\Gamma_{\f O}$ for  orbifolds $\f O$ with $\chi(\f O)>0$ 
whose ramifications collections 
are listed in \eqref{list2}.
Let $\Gamma$ be a finite subgroup  of $Aut(\C\P^1)$. 
Abusing the notation, we will  denote by $\theta_{\Gamma}$  the universal covering $\theta_{\f O}$, where  $\f O$ is an orbifold such that $ \Gamma_{\f O}=\Gamma$.  
Thus, 
$\theta_{\Gamma}$ is defined up to the transformation 
$$\theta_{\Gamma}\rightarrow \delta\circ \theta_{\Gamma},$$ where $\delta$ is a M\"obius transformation.

Let $F$ be a rational function. Recall that two 
decompositions 
of $F$ into compositions of rational functions 
\be \l{ssu}  F=A\circ W\ee 
and 
$$ F=\t A\circ \t W$$ are called  
equivalent if 
\be \l{trans} \t A=A\circ \mu,\ \ \ \ \ \t W=\mu^{-1}\circ W\ee for some M\"obius transformation $\mu$.
Equivalence classes of decompositions of $F$ are in 
a one-to-one correspondence with imprimitivity systems of the monodromy group $Mon(F)$ of $F$.
Namely, if $z_0$ is a non-critical value of $F$ and  $Mon(F)$ is realized as a permutation group acting on the fiber $F^{-1}\{z_0\}$, then to the equivalence class of decomposition \eqref{ssu} corresponds 
the imprimitivity system consisting of 
$d=\deg A$ blocks $W^{-1}\{t_i\},$ $1\leq i \leq d,$ where $\{t_1,t_2,\dots, t_{d}\} =A^{-1}\{z_0\}$ (see e.g. \cite{prime}, Section 2 for more details).

Imprimitivity systems of a transitive permutation group $G$ acting on a set $S$ are in a one-to-one correspondence with subgroups of $G$ containing a stabilizer $G_{\alpha}$ of some element $\alpha$ of $S$ (see 
e.g. \cite{wi}, Theorem 7.5). On the other hand, if $f:\, S\rightarrow \C\P^1$ is a Galois covering, then equality 
\eqref{dega} implies that stabilizer subgroups of $Mon(f)$ are trivial. 
Thus, 
 for any orbifold $\f O$ with $\chi(\f O)>0$ equivalence classes of 
decompositions 
of $\theta_{\f O}: \C\P^1\rightarrow \C\P^1$   are in a one-to-one correspondence with subgroups $\Gamma^{\prime}$ of $\Gamma_{\f O}$, and 
any decomposition of $\theta_{\f O}$ has
the form \be \l{xoo} \theta_{\f O}=A_{\Gamma^{\prime}}\circ  \theta_{\Gamma^{\prime}},\ee
where $\Gamma^{\prime}$  is a
subgroup of the group $\Gamma_{\f O}$ and  $A_{\Gamma^{\prime}}$ is some rational function depending on $\Gamma^{\prime}.$ 
However, the  number  of $\mu$-equivalence classes of compositional 
left factors of $\theta_{\f O}$ in general is less than 
the number of  equivalence classes of 
decompositions of $\theta_{\f O}$. In particular, to conjugate subgroups $\Gamma_1$, $\Gamma_2$ of $\Gamma_{\f O}$ correspond 
$\mu$-equivalent functions $A_{\Gamma_1}$, $A_{\Gamma_2}.$ More precisely, 
the following statement holds.

\bl \l{nml} Let $\f O$ be an orbifold of positive Euler characteristic on $\C\P^1$,  $\Gamma_1$, $\Gamma_2$ subgroups of $\Gamma=\Gamma_{\f O}$, and 
$$\theta_{\f O}=A_{\Gamma_i}\circ  \theta_{\Gamma_i}, \ \ \ i=1,2,$$ corresponding decompositions.
Then  $A_{\Gamma_2}=A_{\Gamma_1}\circ \delta$ for some M\"obius transformation
$\delta$ if and only if $\Gamma_1$, $\Gamma_2$  are conjugate
in $\Gamma$.
\el
\pr 
The conjugacy condition 
$$\Gamma_2=\mu^{-1}\circ \Gamma_1\circ \mu, \ \ \ \mu\in \Gamma,$$ is equivalent to the condition that for any choice of $\theta_{\Gamma_1}$ and $\theta_{\Gamma_2}$ the equality
\be \l{pesik} \delta\circ \theta_{\Gamma_2}= \theta_{\Gamma_1}\circ\mu\ee holds for some M\"obius transformation $\delta.$ 
Assume that \eqref{pesik} holds. Then, since $\theta_{\f O}=\theta_{\f O}\circ \mu$  for any $\mu\in \Gamma$, 
we have: 
\be \l{qaz} \theta_{\f O}=A_{\Gamma_1}\circ \theta_{\Gamma_1}=A_{\Gamma_1}\circ \theta_{\Gamma_1}\circ \mu= A_{\Gamma_1}\circ \delta\circ\theta_{\Gamma_2}.\ee 
Since, on the other hand, \be \l{zaq} \theta_{\f O}=A_{\Gamma_2}\circ \theta_{\Gamma_2},\ee we conclude that 
\be \l{kotik} A_{\Gamma_2}=A_{\Gamma_1}\circ \delta.\ee

In the other direction, assume that \eqref{kotik} holds. Consider the algebraic curve obtained by equating to zero the numerator of $\theta_{\f O}(y)-\theta_{\f O}(x).$
Abusing the notation we will denote this curve simply by 
\be \l{coro} \theta_{\f O}(y)-\theta_{\f O}(x)=0.\ee
Since the rational function $\theta_{\theta}$ is a Galois covering, curve \eqref{coro} 
 splits over $\C(x)$ into a product of factors of degree one
$$y-\mu(x), \ \ \ \mu\in \Gamma.$$ 
On the other hand, since 
\begin{multline*} \theta_{\f O}(y)-\theta_{\f O}(x)=(A_{\Gamma_1}\circ \theta_{\Gamma_1})(y)-(A_{\Gamma_2}\circ \theta_{\Gamma_2})(y)=\\ =(A_{\Gamma_1}\circ \theta_{\Gamma_1})(y)-(A_{\Gamma_1}\circ \delta\circ \theta_{\Gamma_2})(y)
\end{multline*} 
and $y-x$ is a factor of the algebraic curve $A_{\Gamma_1}(y)-A_{\Gamma_1}(x)=0$, 
the algebraic curve 
\be \l{ebs} \theta_{\Gamma_1}(y)-\delta\circ \theta_{\Gamma_2}(x)=0\ee is a factor of curve \eqref{coro}. Therefore, curve \eqref{ebs} also splits over $\C(x)$ into a product of factors of degree one. Since $y$ 
in \eqref{ebs} can be represented locally in the form 
$$y=(\theta_{\Gamma_1}^{-1}\circ \delta\circ \theta_{\Gamma_2})(x),$$
where $\theta_{\Gamma_1}^{-1}$ is a branch of the algebraic function inverse to $\theta_{\Gamma_1},$ we conclude that 
$$\theta_{\Gamma_1}^{-1}\circ \delta\circ \theta_{\Gamma_2}=\mu, \ \ \ \mu\in \Gamma,$$ implying \eqref{pesik}. \qed

\vskip 0.2cm

Lemma \ref{ml} and formula \eqref{vtgh} assure that any rational function $A$ of degree greater than one with  $\chi(\f O_2^A)>0$ is a compositional left factor of some $\theta_{\f O}$ with 
\be \l{hold} \nu(\f O)=\nu(\f O_2^A).\ee Nevertheless, the fact that $A$ is a compositional left factor of  $\theta_{\f O}$ implies only that ${\f O_2^A}\preceq \f O$. 
More precisely, the following statement holds.

\bl \l{lemlas} 
Let $\f O$ be an orbifold of positive Euler characteristic on $\C\P^1$
and $A$ a compositional left factor of  $\theta_{\f O}$ of degree at least two. Then either $\nu(\f O)=\nu(\f O_2^A)$, or one of the following conditions holds:

\begin{itemize}

\item $\nu(\f O)=\{n,n\},$ $n\geq 2,$ and  $\nu(\f O_2^A)=\{d,d\},$  where $d\vert n,$ $d\geq 2.$

\item $\nu(\f O)=\{2,2,n\},$ $n\geq 2$, and $\nu(\f O_2^A)=\{2,2,d\},$  where $d\vert n,$ $d\geq 1.$

\item $\nu(\f O)=\{2,3,3\},$  and   $\nu(\f O_2^A)=\{3,3\}.$ 

\item $\nu(\f O)=\{2,3,4\},$  and  $\nu(\f O_2^A)=\{2,2,3\},$  or  $\nu(\f O_2^A)=\{2,2\}.$

\end{itemize}

\el
\pr 
Since
$\f O_2^{\theta_{\f O}}=\f O,$
it follows from the definition of $\f O_2^A$ and the chain rule that ${\f O_2^A}\preceq \f O$. In particular, 
$\chi({\f O_2^A}) \geq \chi(\f O)>0.$
Since the orbifold $\f O_2^A$ has a universal covering, it cannot have one ramified point, or  two ramified points $z_1,$ $z_2$ such that $\nu(z_1)\neq \nu(z_2).$ Furthermore, it cannot be non-ramified
since any rational function of degree at least two has critical values. 
These observations  imply easily the statements of the lemma. For example, if $\nu(\f O)=\{2,3,3\}$ and $\nu(\f O)\neq \nu(\f O_2^A)$,
then the condition ${\f O_2^A}\preceq \f O$ yields that either $\f O_2^A$ is non-ramified, or $\nu(\f O_2^A)$ is one of the following  collections $\{2\},$ $\{3\},$ $\{2,3\},$ $\{3,3\}.$ 
Thus, excluding orbifolds with no universal covering and the non-ramified sphere, we conclude that $\nu(\f O_2^A)=\{3,3\}.$ 
\qed

\end{section}

\begin{section}{Proof of Theorem \ref{mt1}}

Basing on results of Section 3, in this section we prove Theorem \ref{mt1}. In  more details, for
each orbifold $\f O$ on $\C\P^1$ such that $\chi(\f O)>0$ we list all $\mu$-equivalence classes of rational functions $A$ with $\nu(\f O_2^A)=\nu(\f O)$. We   use the following strategy. First, for each conjugacy class  of subgroups of $\Gamma=\Gamma_{\f O}$ we 
find a compositional left factor $A_{\Gamma^{\prime}}$ of $\theta_{\f O}$
corresponding to a representative $\Gamma^{\prime}$ of this class. Then we reject $A_{\Gamma^{\prime}}$ with  $\nu(\f O_2^{A_{\Gamma^{\prime}}})\neq \nu(\f O)$. Finally, we describe 
$\mu$-equivalence classes of the remaining functions. 
It is clear that any rational function  $A$ of degree one is $\mu$-equivalent to $z^n$ 
for $n=1$ and hence is a Galois covering. So, we will consider only compositional 
left factors  of $\theta_{\f O}$  of degree greater than one. 

\vskip 0.2cm

The following elementary lemma is useful for proving that a concrete rational function $A$ has only three critical values.

\bl \l{laslas} Let $f$ be a rational function of degree $d$ such that  the preimage $f^{-1}\{0,1,\infty\}$ contains at most $d+2$ points. Then $f^{-1}\{0,1,\infty\}$ contains exactly $d+2$ points, 
and $f$ has no critical values distinct from $0,1,$ and $\infty.$
\el 
\pr By the Riemann-Hurwitz formula,
$$2d-2=\sum_{z\in \C\P^1}(\deg_zf-1),$$ implying that 
$$\sum_{z\in f^{-1}\{0,1,\infty\}}(\deg_zf-1)\leq 2d-2,$$
where the equality is attained if and only $f$ has no critical values distinct from $0,1,$ and $\infty.$
Therefore, 
$$\vert f^{-1}\{0,1,\infty\}\vert \geq 
\sum_{z\in f^{-1}\{0,1,\infty\}}\deg_zf- 2d+2=d+2,$$
where the equality is attained if and only $f$ has no critical values distinct from $0,1,$ and $\infty.$ \qed

\begin{subsection}{Functions with $\nu(\f O_2^A)=\{n,n\}$} 

If  $\nu(\f O)=\{n,n\}$, $n\geq 2,$ then 
without loss of generality we may assume that $$\nu(0)=n, \ \ \ \nu(\infty)=n,$$ the group $\Gamma_{\f O}=\Z/n\Z$ is generated by the transformation $$ \alpha: z\rightarrow e^{2\pi i/n} \
z,$$ and   $\theta_{\f O}$ equals
$$\theta_{\Z/n\Z}=z^n,\ \ \ \ \ \ \ n\geq 2.$$
Further, since any subgroup of $\Z/n\Z$ is a cyclic group,  and for any $d\vert n$ the group $\Z/n\Z$ contains
only one cyclic subgroup  of order $n/d$,  
any decomposition of $z^n$ into a composition of rational functions is equivalent to the decomposition \be \l{ega} z^n= z^d \circ  z^{n/d},\ee where $d\vert n$. 
Thus, since  $\f O_2^{z^{d}}=\{d,d\}$  and $\{d,d\}\neq \{n,n\}$
for $d<n$,
we see that 
$\nu(\f O_2^A)=\{n,n\}$ if and only if $A\underset{\mu}{\sim}  z^n.$

\end{subsection} 
\begin{subsection}{Functions with $\nu(\f O_2^A)=\{2,2,n\}$} 
If  $\nu(\f O)=\{2,2,n\}$, $n\geq 2,$ then  we may assume that 
$$\nu(-1)=2, \ \ \ \nu(1)=2, \ \ \ \nu(\infty)=n,$$ 
the group 
$\Gamma_{\f O}=D_n$ is generated by the transformations 
$$\alpha: z\rightarrow e^{2\pi i/n} z, \ \ \  \ \ \ \beta: z\rightarrow \frac{1}{z},$$ and 
$\theta_{\f O}$ equals
\be \l{dn} \theta_{{D_{2n}}}=\frac{1}{2}\left(z^n+\frac{1}{z^n}\right),\ \ \ \ \ \ \ n\geq 2.\ee
Any subgroup $G$ of $D_{2n}$ is either cyclic or dihedral. More precisely, either $G=<\alpha^d>$, where $d\vert n$, or  $G=<\alpha^d,\alpha^r\beta>$, where $d\vert n$ and $0\leq r < d.$  Thus, for any $d\vert n$ there exists one subgroup of the first type and $d$ subgroups of the second. The subgroups of the second type 
form one conjugacy class if $n$ is odd, and one or two
conjugacy classes if $n$ is even according as $d$ is odd or even.

Correspondingly, any decomposition of  \eqref{dn} is equivalent either to the decomposition  
\be \l{ega1}\frac{1}{2}\left(z^n+\frac{1}{z^n}\right)=\frac{1}{2}\left(z^{d}+\frac{1}{z^{d}}\right)\circ z^{n/d},\ee
or to the decomposition \be \l{ega2}
\frac{1}{2}\left(z^{n}+\frac{1}{z^{n}}\right)=\left(\v^{d}T_{d}\right)\circ  \frac{1}{2}\left(\v z^{n/d}+\frac{1}{\v z^{n/d}}\right),
\ee
where $\v^{2d}=1$.
For  $\v$ and $-\v$ 
decompositions \eqref{ega2} are equivalent, since $$T_{d}(-z)=(-1)^{d}T_{d}(z).$$ So, 
for odd $d$ we can assume that 
$\v^{d}=1.$  In either case, 
 $$\v^{d} T_{d}\underset{\mu}{\sim}  T_{d}.$$
 Thus, any compositional left factor $A$ of \eqref{ega2} is $\mu$-equivalent either to $\frac{1}{2}\left(z^{d}+\frac{1}{z^{d}}\right)$ or $T_d$, 
 implying that 
 $\nu(\f O_2^A)=\{2,2,n\}$ if and only if $A$ is $\mu$-equivalent either to II, a) or to II, b).

\end{subsection} 

\begin{subsection}{Functions with $\nu(\f O_2^A)=\{2,3,3\}$} 
Any subgroup of $A_4$ distinct from  $A_4$ is isomorphic to one of the following groups:  $\{e\},$ $\Z/2\Z,$ $\Z/3\Z$, $D_4$. Furthermore, any two isomorphic subgroups are conjugate in $A_4.$ 
Thus, the function \be \l{a4} \theta_{A_4}= -\frac{1}{64}\frac{z^3(z^3-8)^3}{(z^3+1)^3},\ee 
which is a universal covering of the orbifold $\f O$ defined by the equalities
\be \l{barsu} \nu(0)=3, \ \ \ \ \nu(1)=2, \ \ \ \ \nu(\infty)=3,\ee
has, up to the change 
$A\rightarrow A\circ \mu,$ where $\mu\in \Aut(\C\P^1)$,  three compositional left factors of degrees 6, 4, and 3, correspondingly.  Moreover, these factors cannot be $\mu$-equivalent since 
they have different degrees. 

Considering 
the obvious decomposition 
$$-\frac{1}{64}\frac{z^3(z^3-8)^3}{(z^3+1)^3}=-\frac{1}{64}\frac{z(z-8)^3}{(z+1)^3}\circ z^3$$
and the decomposition
\be \l{dee} -\frac{1}{64}\frac{z^3(z^3-8)^3}{(z^3+1)^3}=-\frac{1}{64}z^3\circ \frac{z^2-4}{z-1}\circ \frac{z^2+2}{z+1}\ee
found in \cite{mp2}, we see that these factors  are 
\be \l{fa4} -\frac{1}{64}\left(\frac{z^2-4}{z-1}\right)^3,\ \ \ \ \ \  -\frac{1}{64}\frac{z(z-8)^3}{(z+1)^3}, \ \ \ \ \ \   -\frac{1}{64}z^3. \ee
Since Lemma \ref{lemlas} implies that a
 compositional left factor $A$ of  $\theta_{\f O}$ 
satisfies $\nu(\f O)=\nu(\f O_2^A)$ unless $A\sim z^3,$  we conclude that
  a rational function $A$ satisfies  $\nu(\f O_2^A)=\{2,3,3\}$ if and only if $A$
is $\mu$-equivalent to one of the functions listed in III).

Notice that  function \eqref{a4} does not coincide with the function 
\be \l{enott} \left(\frac{z^4+2i\sqrt{3}z^2+1}{z^4-2i\sqrt{3}z^2+1}\right)^3\ee
found by Klein. However, it is easy to see that function \eqref{a4} 
along with \eqref{enott} is a universal covering of $\f O$ given by \eqref{barsu}. 
Indeed, by the uniqueness of the universal covering, 
it is enough to show that  \eqref{a4} satisfies the following conditions:
$f$ has only three critical values $0,1$ and $\infty$, the multiplicity of any  critical point of $f$ over $0$ or 
$\infty$ is 3,  and the multiplicity of any critical point of $f$ over $1$ is 2. Clearly, 
equalities \eqref{a4} 
and  \be \l{egik} f-1=-{\frac {1}{64}}\,{\frac { \left( {z}^{6}+20\,{z}^{3}-8 \right) ^{2}}{
 \left( {z}^{3}+1 \right) ^{3}}}
\ee imply that $f^{-1}\{0,1,\infty\}$ contains at most 14 points, implying by Lemma \ref{laslas} that $f^{-1}\{0,1,\infty\}$ contains exactly 14 points,
and $f$ has no critical values distinct from  $0,1,$ and $\infty$.
It follows now from \eqref{a4} 
and \eqref{egik} that  $f$ has the required ramification  over $0$, $1$, and $\infty$.

Although what is written above  certainly {\it proves} that any rational function with $\nu(\f O_2^A)=\{2,3,3\}$ has one of the forms listed in III), the reader may ask how to {\it find} function  \eqref{a4} and its compositional left factors. 
A  convenient framework  for this purpose is provided by the ``Dessins d'enfants'' theory which 
interprets functions $\theta_{\f O}$ for orbifolds $\f O$ with $\chi(\f O)>0$ as Belyi functions of Platonic solids (see \cite{couv}, \cite{mazv}).
Assuming that the reader is familiar with rudiments of this theory (see e.g. the books \cite{lazv}, \cite{des1}), below we  sketch the 
corresponding calcu\-lations. 

It follows from $\Z/3\Z< A_4$  that  a universal covering of the orbifold given by \eqref{barsu} can be written in the form
 \be \l{ox} \theta_{\f O}=A\circ \theta_{\Z/3\Z}=A\circ z^3\ee 
for some rational function $A$, and  the chain rule implies  that $A$
is a Belyi function. The dessin $\lambda$ corresponding to the Belyi function $ \theta_{\f O}$ is 
the  tetrahedron shown on Fig. \ref{a4_1}, 
\begin{figure}[ht]
\medskip
\epsfxsize=2.3truecm
\centerline{\hskip -0.2cm\epsffile{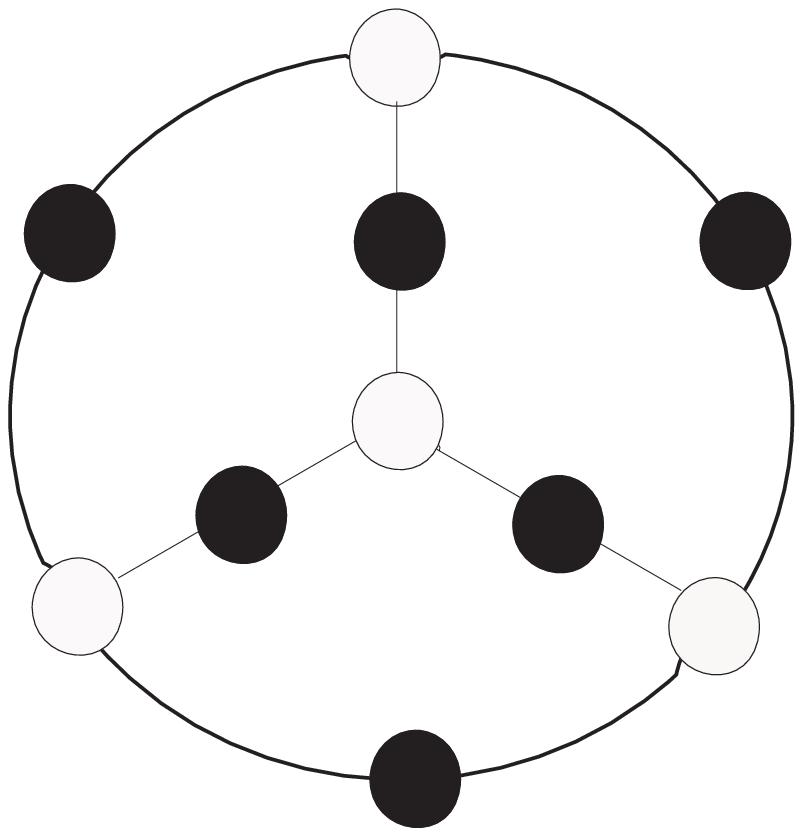}}
\caption{}\l{a4_1}
\medskip
\end{figure}
\noindent
where  as usual  white vertices are preimages of 0, black vertices are preimages of 1, and ``centers'' of faces are preimages of $\infty$.
If the ``interior'' white vertex of $\lambda$
is placed at the origin while the center of the ``exterior'' face is placed at infinity, then  the dessin  
$\lambda_A$ corresponding to the Belyi function $A$ is obtained 
from $\lambda$
by factoring  through the action of the group $\Z/3\Z$ viewed as a rotation group of order three around the origin.
The corresponding dessin $\lambda_A$ is shown on Fig. \ref{a4_2}. 

\begin{figure}[ht]
\medskip
\epsfxsize=2.4truecm
\centerline{\hskip -0.2cm\epsffile{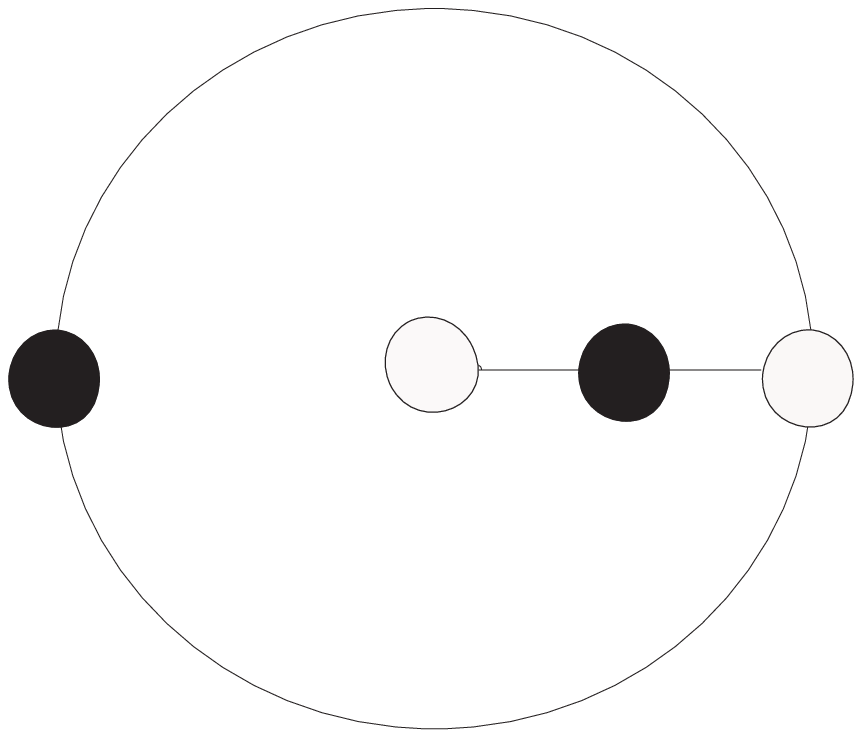}}
\caption{}\l{a4_2}
\medskip
\end{figure}
\noindent By construction, the white vertex of valency 1 of $\lambda_A$ is 0 and  the center of the face of valency 1 is $\infty.$
However, we still can place the center of the interior face of $\lambda_A$ arbitrary, say at the point -1. Then 
\be \l{rfd} A=\frac{az(z-b)^3}{(z+1)^3}\ee for some $a,b\in \C.$
Finally, since the finite roots of the derivative of \eqref{rfd} are $b$ and the roots $\alpha_1,\alpha_2$ of the polynomial $z^2+(2b+4)z-b,$ it follows  from the  
conditions $A(\alpha_1)=A(\alpha_2)=1$ and $\alpha_1\neq \alpha_2$ that $a=-1/64$ and
$b=8$. Thus, we arrive to formula \eqref{a4} and the second function in \eqref{fa4}.

Similarly, the inclusion $\Z/2\Z< A_4$ implies that a universal covering of the orbifold given by \eqref{barsu} can be written in the form
\be \l{ox1} \theta_{\f O}=B\circ z^2\ee 
for some Belyi function $B$, as in formula \eqref{enott}. However, 
in order to view the automorphism of order two of the dessin shown on Fig. \ref{a4_1} 
as a rotation of the second order about the origin, we must redraw it 
placing one of its black vertices at infinity as it is shown on Fig. \ref{a4+}. 
\begin{figure}[ht]
\medskip
\epsfxsize=5.3truecm
\centerline{\hskip 0.3cm\epsffile{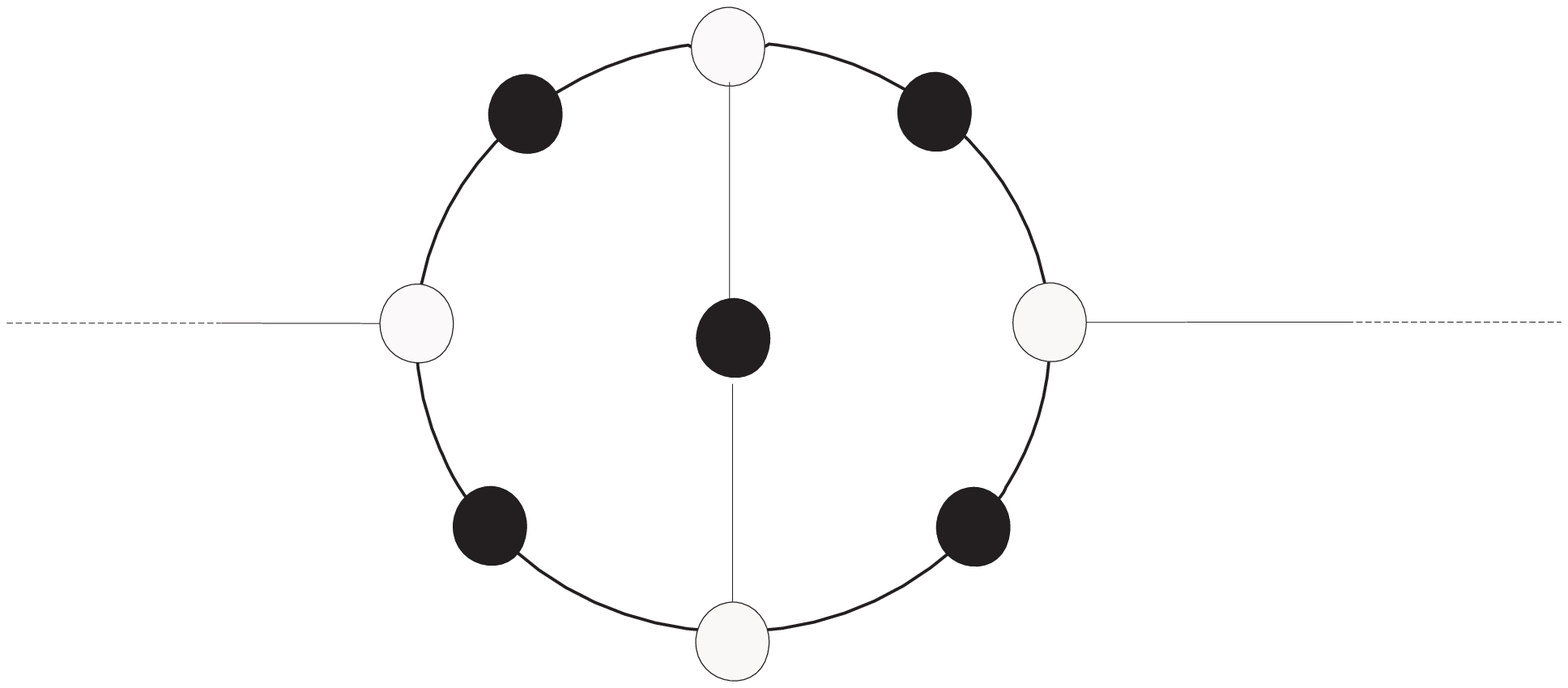}}
\caption{}\l{a4+}
\medskip
\end{figure}
Factoring now through $\Z/2\Z$, we see that  the dessin $\lambda_B$ corresponding to $B$ 
is the one depicted on Fig. \ref{a4_3}. 
\begin{figure}[ht]
\medskip
\epsfxsize=3.1truecm
\centerline{\hskip -1cm\epsffile{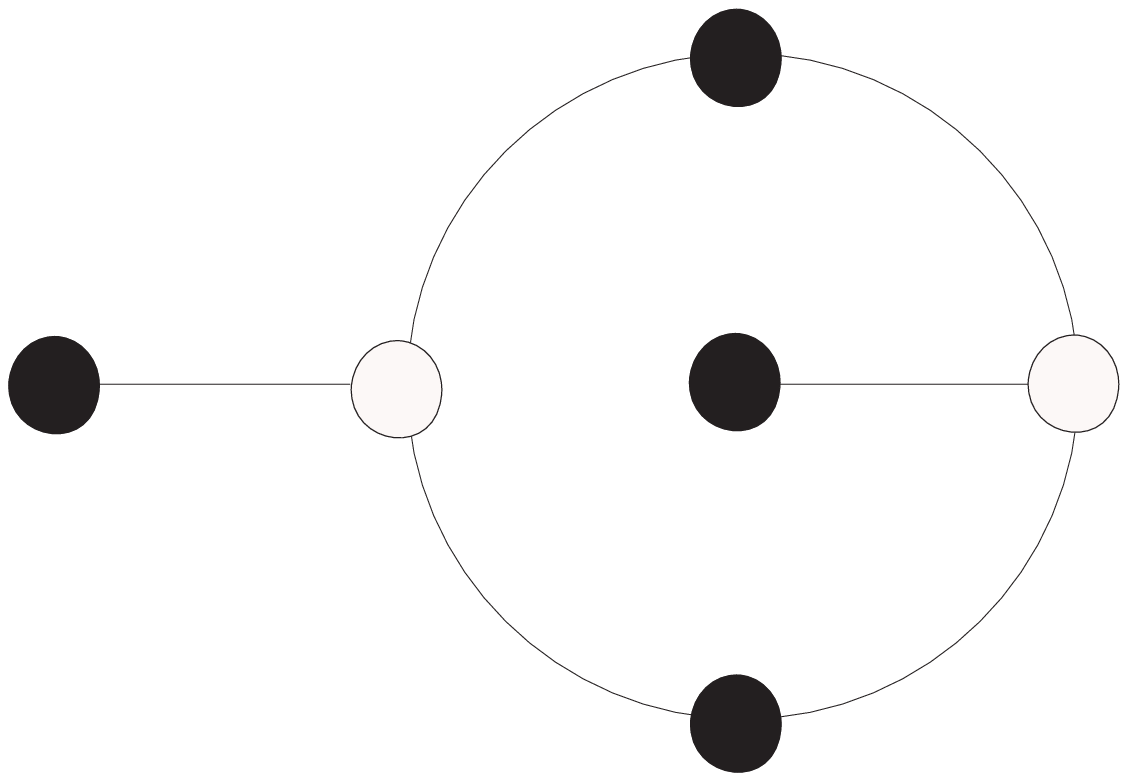}}
\caption{}\l{a4_3}
\medskip
\end{figure}

We can find the function $B$ using a reasoning similar to the one used for finding $A$. However, 
one can reduce calculations 
using the fact that 
the decomposition 
\be \l{xro} \theta_{A_4}=-\frac{1}{64}z^3\circ \frac{z(z^3-8)}{(z^3+1)}\ee  corresponds to  the 
 subgroup  $D_4$. Since any subgroup of $ A_4$ isomorphic to $\Z/2\Z$ is contained in $D_4$, 
this implies that $B\underset{\mu}{\sim}\t B$, where 
 $$\t B=-\frac{1}{64}z^3\circ f$$ for some rational function $f$ of degree two.
Since we can place the centers of the faces of $\lambda_B$ at 1 and $\infty$,
and assume that the sum of two white vertices of valency 3 of $\lambda_B$ is zero, we have:
$$\t B=-\frac{1}{64}\left(\frac{z^2-c}{z-1}\right)^3, \ \ \ c\in \C.$$ Finally, since the finite roots of $\t B^{\prime}(z)$ are $\pm \sqrt{c}$ and $1\pm \sqrt{1-c}$, it follows from the condition $$\t B(1+ \sqrt{1-c})=\t B(1-\sqrt{1-c})=1$$
that $c=4$. 
This gives us the first function in \eqref{fa4}.

\end{subsection} 

\begin{subsection}{Functions with $\nu(\f O_2^A)=\{2,3,4 \}$} Any subgroup of  $S_4$ distinct from  $S_4$  is isomorphic to one of the following groups: 
$\{e\},$ $\Z/2\Z,$ $\Z/3\Z$, $\Z/4\Z$, $D_4$, $S_3,$ $D_8$, $A_4$. Furthermore, $S_4$ has {\it two} conjugacy classes of subgroups isomorphic to $\Z/2\Z$, and {\it two} conjugacy classes of  subgroups isomorphic to $D_4$. Thus,  a universal covering $\theta_{\f O}$ of an 
 orbifold $\f O$ with $\nu(\f O)=\{2,3,4 \}$
has, up to the change 
$A\rightarrow A\circ \mu,$ where $\mu\in \Aut(\C\P^1)$,  ten compositional left factors.  Nevertheless, not all of them satisfy  $\nu(\f O_2^A)=\{2,3,4\}$. For example, for the factors $A$ of degree two and three, 
corresponding to the subgroups $A_4$ and $D_8$, clearly 
$\nu(\f O_2^A)\neq \{2,3,4\}$ since a rational function of degree less than 4  cannot have a critical point of multiplicity 4. Besides, 
we will show that one of the factors  corresponding to $D_4$  has the ramification orbifold $\{2,2,3\}$.
As above, compositional left factors of $\theta_{S_4}$ can be found with the use of ``Dessins d'enfants'' theory. 
We  omit the details of calculations restricting ourselves  by a
formal proof of Theorem \ref{mt1}.

First, it was established already by Klein that the function 
\be \l{xorr}  \theta_{S_4}=\frac{1}{108}\frac{(x^8+14x^4+1)^3}{x^4(x^4-1)^4}\ee
is a universal covering of the orbifold $\f O$ defined by the equalities
$$ \nu(0)=3, \ \ \ \ \nu(1)=2, \ \ \ \ \nu(\infty)=4.$$ 
Clearly, the compositional left factor 
in the decomposition 
$$ \theta_{S_4}=
\frac{1}{54}\frac{(z+7)^3}{(z-1)^2}\circ \theta_{D_8},$$  where 
$\theta_{D_8}$ is given by \eqref{dn}, 
 corresponds  to the subgroup $D_8$. 
Using now the decomposition $$\theta_{D_8}=\frac{1}{2}\left(z+\frac{1}{z}\right)\circ z^4,$$ we obtain the compositional left factor 
 \be \l{xlop} R=
L\circ \frac{1}{2}\left(z+\frac{1}{z}\right)={\frac {1}{108}}\,{\frac { \left( {z}^{2}+14\,z+1 \right) ^{3}}{z
 \left( z-1 \right) ^{4}}}.
\ee of $ \theta_{S_4}$ corresponding 
to   $\Z/4\Z$. Since $R$ has a pole of order four, it follows from Lemma \ref{lemlas}  that $\nu(\f O_2^{R})=\{2,3,4\}.$

Similarly, using the decompositions $$\frac{1}{2}\left(z^4+\frac{1}{z^4}\right)=T_2\circ \frac{1}{2}\left(z^2+\frac{1}{z^2}\right),$$ $$\frac{1}{2}\left(z^4+\frac{1}{z^4}\right)=-T_2\circ \frac{1}{2}\left(iz^2+\frac{1}{iz^2}\right),$$
we obtain the  compositional left factors 
\be \l{ep} B_1=L\circ T_2=\frac{1}{27}\,{\frac { \left( {z}^{2}+3 \right) ^{3}}{ \left( z^2-1 \right) ^{2}
}},\ee
\be \l{epp}
B_2=L\circ (-T_2)=-\frac{1}{27}\,{\frac { \left( z^2-4 \right) ^{3}}{{z}^{
4}}}
\ee of $ \theta_{S_4}$, 
corresponding to  subgroups  isomorphic to $D_4$. Since $B_2$ has a pole of order four, $\nu(\f O_2^{B_2})=\{2,3,4\}.$
 On the other hand, it follows from the equalities \eqref{ep} and 
$$
B_1-1=\frac{1}{27}\,{\frac {{z}^{2} \left( z^2-9 \right) ^{2}}
{ \left( z^2-1 \right) ^{2} }},
$$ that 
$\nu(\f O_2^{B_1})=\{2,2,3\}$. The subgroups  corresponding to the functions $B_1$ and $B_2$  
are not conjugated since otherwise Lemma \ref{nml} would imply that  $ B_1= B_2\circ \mu$ for some $\mu\in Aut(\C\P^1)$ in contradiction with 
$\nu(\f O_2^{B_1})\neq \nu(\f O_2^{B_2})$.

Further, composing  $L$ with compositional left factors $T_4$, $-T_4$ of $\theta_{D_8}$
we obtain the  compositional left factors
 \be \l{prev}  B_3=L\circ T_4= 
{\frac {4}{27}}\,{\frac { \left( {z}^{4}-{z}^{2}+1 \right) ^{3}}{{z}^{
4} \left( z^2-1 \right) ^{2}}},
\ee 
\be \l{prevv}  B_4=L\circ (- T_4)= -\frac{1}{27}\,{\frac { \left( 2\,{z}^{2}+1 \right) ^{3} \left( 2\,{z}^{2}-3
 \right) ^{3}}{ \left( 2\,{z}^{2}-1 \right) ^{4}}}
\ee 
of $ \theta_{S_4}$, corresponding to subgroups isomorphic to   $\Z/2\Z$. Moreover, the subgroups  corresponding to $B_3$ and $B_4$  
are not conjugated, since the function $B_3$ has the ramification $2,2,4,4$
 over infinity,
while $B_4$ has the ramification $4,4,4$.

The function 
$$M=\\ -{\frac {256}{27}}\,{z}^{3} \left( z-1 \right)$$ corresponding to the subgroup $ S_3\cong D_6$ is obtained from the decomposition  
\begin{multline} \l{esz} 
-{\frac {1}{432}}\,{\frac { \left( 16\,{z}^{8}-56\,{z}^{4}+1 \right) ^
{3}}{{z}^{4} \left( 4\,{z}^{4}+1 \right) ^{4}}}= \\ -{\frac {256}{27}}\,{z}^{3} \left( z-1 \right) \circ \frac{1}{8}\,{\frac { \left( 2\,{z}^{2}+2\,z-1 \right)  \left( 4\,{z}^{4}+8\,{
z}^{2}+1 \right) }{z \left( 4\,{z}^{4}+1 \right) }}.
\end{multline}
The function in the left part of \eqref{esz} is obtained from \eqref{xorr} by the substitution $z=\omega z,$ where $\omega^4=-4$.

Finally, consider the decomposition \be \l{bs4}  
256\frac{{z}^{3} \left( {z}^{6}-7z^3-8 \right) ^{3}}
{\left( z^6 +20z^3-8\right) ^{4}}= -\frac {4x}{{x}^{2}+1-2\,x}\circ \theta_{A_4},\ee
where 
$\theta_{A_4}$ is given by formula \eqref{a4}.  
The function $f$ in the left part of equality \eqref{bs4} is
$\mu$-equivalent to \eqref{xorr}. This fact can be checked as above using the 
formula
$$ f-1=-{\frac { \left( {z}^{2}+2 \right) ^{2} \left( {z}^{4}-2\,{z}^{2}+4
 \right) ^{2} \left( {z}^{2}-4\,z-2 \right) ^{2} \left( {z}^{4}+4\,{z}
^{3}+18\,{z}^{2}-8\,z+4 \right) ^{2}}{\left( z^6 +20z^3-8\right) ^{4} }}
$$ and Lemma \ref{laslas}. Composing now the left factor of $f$ from decomposition \eqref{bs4} with 
the  left factor of  $\theta_{A_4}$ of degree 4 found above, we arrive to the function 
$$\left(-\frac {4x}{{x}^{2}+1-2\,x}\right)\circ\left( -\frac{1}{64}\frac{z(z-8)^3}{(z+1)^3}\right)=256\,{\frac {z \left( z^2-7z-8 \right) ^{3}}{
 \left( {z}^{2}+20\,z-8 \right) ^{4}}}
$$  corresponding to the subgroup $\Z/3\Z$.

\end{subsection}

\begin{subsection}{Functions with $\nu(\f O_2^A)=\{2,3,5\}$} 
The subgroups of $A_5$ distinct from $A_5$  are $\{e\}$, $\Z/2\Z,$ $\Z/3\Z$, $D_4$, $\Z/5\Z$, $D_6$, $D_{10}$, and $A_4$. Since any two isomorphic subgroups in $A_5$ are conjugate, 
it follows from Lemma \ref{nml} that 
 a universal covering $\theta_{\f O}$ of an 
 orbifold $\f O$ with $\nu(\f O)=\{2,3,5 \}$ has, 
up to the transformation 
$A\rightarrow A\circ \mu,$ where 
$\mu\in \Aut(\C\P^1)$, 
eight compositional left factors $A$  of degrees 60, 30, 20, 15, 12, 10, 6, and 5, correspondingly. Since these factors have different degrees, they cannot be $\mu$-equivalent.
Furthermore, by Lemma \ref{lemlas}, all these factors satisfy the condition $\nu(\f O_2^A)=\{2,3,5\}$. 
Therefore,  
 up to the $\mu$-equivalence, there exist exactly eight rational  functions $A$ with  $\nu(\f O_2^A)=\{2,3,5\}$, and  in order to finish the proof of Theorem \ref{mt1} we only must check that 
all functions $A$ from list V) satisfy the condition $\nu(\f O_2^A)=\{2,3,5\}$. In  turn, the last statement follows easily from Lemma \ref{laslas} and formulas for $A-1$ given below.

$${\frac {1}{1728}}\,{\frac { \left( {z}^{30}-522\,{z}^{25}-10005\,{z}^{
20}-10005\,{z}^{10}+522\,{z}^{5}+1 \right) ^{2}}{{z}^{5} \left( {z}^{
10}-11\,{z}^{5}-1 \right) ^{5}}}
\leqno{a)}$$

$$-{\frac {1}{6144}}\, \left( 3\,z+11 \right)  \left( 3\,{z}^{2}+2\,z+27
 \right) ^{2}\leqno{b)}$$

$${\frac {1}{1728}}\,{\frac { \left( {z}^{2}+12\,z+40 \right)  \left( {z
}^{2}-6\,z+4 \right) ^{2}}{z-5}}\leqno{c)}$$

$$
-{\frac {1}{9}}\,{\frac { \left( 180\,{z}^{2}+380\,z+229 \right)  \left( 20\,{z}^
{2}+20\,z+41 \right) ^{2} \left( 20\,{z}^{2}-580\,z-979 \right) ^{2}}{
 \left( 20\,{z}^{2}+140\,z+101 \right) ^{5}}}
\leqno{d)}$$

$${\frac {1}{1728}}\,{\frac { \left( {z}^{6}-522\,{z}^{5}-10005\,{z}^{
4}-10005\,{z}^{2}+522\,{z}+1 \right) ^{2}}{{z} \left( {z}^{
2}-11\,{z}-1 \right) ^{5}}}
\leqno{e)}$$

$$
\scriptstyle
-{\frac {1}{27}}\,{\frac { \left( 10\,z+3 \right)  \left( 20\,{z}^{2}+20\,z+1
 \right)  \left( 10\,{z}^{2}+10\,z+3 \right) ^{2} \left( 500\,{z}^{4}+
300\,{z}^{3}+70\,{z}^{2}+10\,z+1 \right) ^{2}}{ \left( 20\,{z}^{2}+10
\,z+1 \right) ^{5}}}
\leqno{f)}$$

$$\scriptstyle
-{\frac {1}{64}}\,{\frac { \left( {z}^{2}+5 \right) ^{2} \left( 8\,{z}
^{4}-100\,{z}^{3}+2055\,{z}^{2}+500\,z+200 \right) ^{2} \left( {z}^{4}
-350\,{z}^{3}-2190\,{z}^{2}+1750\,z+25 \right) ^{2}}{ \left( {z}^{4}+
55\,{z}^{3}-165\,{z}^{2}-275\,z+25 \right) ^{5}}}
\leqno{g)}$$

$$\scriptstyle
{\frac {1}{1728}}\,{\frac { \left( {z}^{2}+4 \right)  \left( {z}^{2}-2
\,z-4 \right) ^{2} \left( {z}^{4}+3\,{z}^{2}+1 \right) ^{2} \left( {z}
^{4}+6\,{z}^{3}+21\,{z}^{2}+36\,z+61 \right) ^{2} \left( {z}^{4}-4\,{z
}^{3}+21\,{z}^{2}-34\,z+41 \right) ^{2}}{ \left( z-1 \right) ^{5}
 \left( {z}^{4}+{z}^{3}+6\,{z}^{2}+6\,z+11 \right) ^{5}}}
\leqno{h)}$$

\end{subsection}

\end{section}

\begin{section}{Functions with $\chi(\f O_2^A)=0$}
Let $\f O$ be an orbifold on $\C\P^1$ such that $\chi(\f O)= 0$. Then 
the corresponding group $\Gamma_{\f O}$
is generated by translations of $\C$ by elements of some lattice $L\subset \C$ of rank two and the transformation $z\rightarrow  \v z,$ where $\v$ is an $n$th root of unity with $n$ equal to 2,3,4, or 6, such that  $\v L=L.$ 
We will denote  by $\Lambda_{\f O}$ the subgroup of $\Gamma_{\f O}$  generated by translations. 
The group $\Lambda_{\f O}$ is normal in $\Gamma_{\f O}$, and can be described as the kernel of the homomorphism $\psi: \Gamma_{\f O} \rightarrow \C$ which sends $\sigma=az+b\in \Gamma_{\f O}$ 
 to $\psi(\sigma)=a\in \C.$
For  $\nu(\f O)=\{2,2,2,2\}$ the complex structure of $\C/\Lambda_{\f O}$ may be arbitrary, and the function  $\theta_{\f O}$ is the corresponding Weierstrass function $\wp(z).$ On the other hand, for $\nu(\f O)$ equal $\{2,4,4\}$, $\{3,3,3\}$, or $\{2,3,6\}$  the complex structure of $\C/\Lambda_{\f O}$  is   
rigid and arises from the tiling of $\C$ by squares, equilateral triangles, or alternately colored equilateral triangles, respectively. Accordingly, the function $\theta_{\f O}$ 
may be written in terms of the corresponding
Weierstrass functions as $\wp^2(z),$ $\wp^{\prime }(z),$ and $\wp^{\prime 2}(z)$ (see  \cite{mil2} and \cite{fk}, Section IV.9.12).

The following statement provides a geometric description of covering maps $A: \f  O_1\rightarrow \f O_2$  between  orbifolds of zero characteristic.

\bt \l{mt3} Let $A$ be a rational functions. Then $A:\f O_1\rightarrow \f O_2$ is a covering map 
between some orbifolds of zero Euler characteristic $\f O_1$ and $\f O_2$ on $\C\P^1$ if and only if 
there exist elliptic curves
$\f C_1$ and $\f C_2$, subgroups  $\Omega_1\subseteq {\rm Aut}({\f  C_1})$ and  $\Omega_2\subseteq {\rm Aut}({\f  C_2})$, and a holomorphic map $\alpha: \f C_1\rightarrow \f C_2$
such that the diagram 
\be \l{xxuu}
\begin{CD}
\f C_1 @>\alpha>> \f C_2 \\
@VV\pi_1 V @VV\pi_2 V\\ 
\C\P^1 @>A >> \ \ \C \P^1\, ,
\end{CD}
\ee
where $\pi_1: {\f C_1}\rightarrow   
{\f C_1}/\Omega_1$ and $\pi_2: {\f C_2}\rightarrow   
{\f C_2}/\Omega_2$ are quotient maps, is commutative. 
\et
\pr
If $A$ is a rational function such $A:\f O_1\rightarrow \f O_2$ is a covering map between some 
orbifolds of zero Euler characteristic $\f O_1$ and $\f O_2$,
then there exists an isomorphism $F=az+b$, $a,b\in \C,$ of the complex plane which makes the diagram
\be \l{diakot}
\begin{CD}
\C @>F=az+b>> \C \\
@VV\theta_{\f O_1}V @VV\theta_{\f O_2}V\\ 
\f O_1 @>A >> \f O_2\ 
\end{CD}
\ee
commutative and 
satisfies \eqref{homm}
for some homomorphism $\phi:\, \Gamma_{\f O_1}\rightarrow \Gamma_{\f O_2}$. 
Moreover, $\phi$ is a monomorphism since $F$ is invertible and 
hence the 
equality $F\circ \sigma =F$ implies that $\sigma=z.$

It is clear that $\f C_1=\C/\Lambda_{\f O_1}$ and $\f C_2=\C/\Lambda_{\f O_2}$ are Riemann surfaces of genus
one, and the groups 
$$\Omega_1\cong \Gamma_{\f O_1}/\Lambda_{\f O_1}, \ \ \ \ \Omega_2\cong \Gamma_{\f O_2}/\Lambda_{\f O_2}$$ 
are cyclic groups of order 2,3,4, or 6. Moreover, 
we can consider $\f C_1$ and $\f C_2$ as elliptic curves, whose marked points are projections of the origin, and 
$\Omega_1$ and  $\Omega_2$ as automorphism groups of $\f C_1$ and $\f C_2$.
Further, 
$$\theta_{\f O_1}=\pi_1\circ \psi_1, \ \ \ \ \theta_{\f O_2}=\pi_2\circ \psi_2,$$ where $$\psi_1:\C\rightarrow \C/\Lambda_{\f O_1}\cong \f C_1, \ \ \ \psi_2:\C\rightarrow \C/\Lambda_{\f O_2}\cong \f C_2$$ 
and 
$$\pi_1:\f C_1\rightarrow \f C_1/\Omega_1\cong \C\P^1, \ \ \ \ \pi_2:\f C_2\rightarrow \f C_2/\Omega_2\cong \C\P^1,$$ 
are quotient maps. Finally, 
since $\phi$ is a monomorphism, it maps elements of infinite order of $\Gamma_{\f O_1}$ to elements of infinite order of $\Gamma_{\f O_2}$.  
Thus, $\phi(\Lambda_{\f O_1})\subset \Lambda_{\f O_2},$ implying that $F$ descends to a holomorphic map $\alpha: \f C_1\rightarrow \f C_2$ which makes the 
diagram 
\be 
\begin{CD} \l{gpa3}
\C @>F=ax+b>> \C \\
@VV  \psi_1 V @VV  \psi_2  V\\ 
\f C_1 @>\alpha>> \f C_2\\
@VV \pi_1 V @VV \pi_2 V\\ 
\C\P^1 @>A >> \ \ \C\P^1\,\  
\end{CD}
\ee
commutative. 

In the other direction, we can complete any diagram \eqref{xxuu} to diagram \eqref{gpa3}, setting $\psi_1$ and $\psi_2$ equal the usual universal coverings of the Riemann surfaces $\f C_1$ and $\f C_2$. 
Since 
 $\pi_1$ and $\pi_2$ are Galois coverings
and the maps $\psi_1$ and $\psi_2$  are non-ramified, it is easy to see that the maps 
$\pi_1\circ \psi_1:\, \C \rightarrow \f O_2^{\pi_1}$ and $\pi_2\circ \psi_2:\, \C \rightarrow \f O_2^{\pi_2}$ are universal coverings of the  orbifolds $\f O_2^{\pi_1}$ and $\f O_2^{\pi_2}$, implying that  
 $A:\f O_2^{\pi_1}\rightarrow \f O_2^{\pi_2}$ is a covering map between orbifolds. Finally, since 
$$\t{\f O_2^{\pi_1}}= \t{\f O_2^{\pi_2}}=\C,$$ 
 these orbifolds have zero Euler characteristic.  
 \qed

\vskip 0.2cm

Obviously, Corollary \ref{cor2} and Theorem \ref{mt3} imply the first part of Theo\-rem \ref{mt2}.  On the other hand, in order to prove the second part we must show that if  $A:\f O_1\rightarrow \f O_2$ is 
a covering map between orbifolds of zero Euler characteristic  such that $\chi(\f O_2^A)>0,$ then $A$
is $\mu$-equivalent either to a cyclic  function for some $n\leq 4$, or to a dihedral function for some  $n\leq 4$, or
to a tetrahedral function.  The theorem below provides  a more precise version of the required statement.

\bt \l{l}  Let $A$ be a rational function and $\f O_1$, $\f O_2$ orbifolds such that $\chi(\f O_1)=\chi(\f O_2)=0$.  
Assume that  
$A:\, \f O_1\rightarrow \f O_2$ is a covering map between orbifolds. Then either $\chi(\f O_2^A)=0$ and 
$\f O_2=\f O_2^A$, $\f O_1=\f O_1^A$, or $\chi(\f O_2^A)>0$ and one of the following conditions holds:

$$\nu(\f O_2^A)=\{2,2\}, \ \   {A}\underset{\mu}{\sim}z^2, \ \ \nu(\f O_1)=\nu(\f O_2)=\{2,2,2,2\},\leqno{1.}$$
$$\nu(\f O_2^A)=\{2,2\}, \ \   {A}\underset{\mu}{\sim}z^2,  \ \ \nu(\f O_1)=\nu(\f O_2)=\{2,4,4\},\leqno{2.}$$ 
$$\nu(\f O_2^A)=\{2,2\}, \ \   {A}\underset{\mu}{\sim}z^2,  \ \ \nu(\f O_1)=\{2,2,2,2\},\ \ \nu(\f O_2)=\{2,4,4\},\leqno{3.}$$ 
$$\nu(\f O_2^A)=\{2,2\}, \ \   {A}\underset{\mu}{\sim}z^2,  \ \ \nu(\f O_1)=\{3,3,3\}, \ \ \nu(\f O_2)=\{2,3,6\},\leqno{4.}$$ 
$$\nu(\f O_2^A)=\{3,3\}, \ \   {A}\underset{\mu}{\sim}z^3,  \ \ \nu(\f O_1)=\nu(\f O_2)=\{3,3,3\},\leqno{5.}$$  
$$\nu(\f O_2^A)=\{3,3\}, \ \   {A}\underset{\mu}{\sim}z^3,  \ \ \nu(\f O_1)=\{2,2,2,2\}, \ \ \nu(\f O_2)=\{2,3,6\},\leqno{6.}$$ 
$$\nu(\f O_2^A)=\{4,4\}, \ \   {A}\underset{\mu}{\sim}z^4,  \ \ \nu(\f O_1)=\{2,2,2,2\},\ \ \nu(\f O_2)=\{2,4,4\},\leqno{7.}$$

$$\nu(\f O_2^A)=\{2,2,2\}, \ \   {A}\underset{\mu}{\sim} \frac{1}{2}\left(z^2+\frac{1}{z^2}\right),\ \    \nu(\f O_1)=\nu(\f O_2)=\{2,2,2,2\},\leqno{8.}$$

$$\hskip -0.25cm \nu(\f O_2^A)=\{2,2,2\}, \,  {A}\underset{\mu}{\sim} \frac{1}{2}\left(z^2+\frac{1}{z^2}\right),\, \nu(\f O_1)=\{2,2,2,2\},\,  \nu(\f O_2)=\{2,4,4\},\leqno{9.}$$

$$\hskip -0.3cm\nu(\f O_2^A)=\{2,2,3\}, \    {A}\underset{\mu}{\sim} \frac{1}{2}\left(z^3+\frac{1}{z^3}\right), \   \nu(\f O_1)=\{3,3,3\},\  \nu(\f O_2)=\{2,3,6\},\leqno{10.}$$ 
$$\nu(\f O_2^A)=\{2,2,3\},  \ \ \ \   {A}\underset{\mu}{\sim} T_3,  \ \ \ \  \nu(\f O_1)=\nu(\f O_2)=\{2,3,6\},\leqno{11.}$$ 
$$\hskip -0.25cm  \nu(\f O_2^A)=\{2,2,4\},\,  {A}\underset{\mu}{\sim} \frac{1}{2}\left(z^4+\frac{1}{z^4}\right),\,  \nu(\f O_1)=\{2,2,2,2\},\, \nu(\f O_2)=\{2,4,4\},\leqno{12.} $$ 
$$\hskip -0.3cm  \nu(\f O_2^A)=\{2,2,4\},\ \ \ {A}\underset{\mu}{\sim} T_4,\ \ \ \nu(\f O_1)=\{2,2,2,2\},\ \ \ \nu(\f O_2)=\{2,4,4\},\leqno{13.} $$ 
$$\nu(\f O_2^A)=\{2,2,4\}, \ \ \  {A}\underset{\mu}{\sim} T_4,\ \ \ \nu(\f O_1)=\nu(\f O_2)=\{2,4,4\},\leqno{14.}$$ 
$$\nu(\f O_2^A)=\{2,3,3\},\ \ \ {A}\underset{\mu}{\sim} -\frac{1}{64}\frac{z(z-8)^3}{(z+1)^3},\ \ \ \nu(\f O_1)=\nu(\f O_2)=\{2,3,6\}.\leqno{15.}$$ 
$$ \nu(\f O_2^A)=\{2,3,3\}, \ \ \ \  {A}\underset{\mu}{\sim} -\frac{1}{64}\left(\frac{z^2-4}{z-1}\right)^3,\leqno{16.}$$  
$$\nu(\f O_1)=\{2,2,2,2\}, \ \ \ \ \nu(\f O_2)=\{2,3,6\},$$  
$$ \nu(\f O_2^A)=\{2,3,3\}, \ \ \ \   A\underset{\mu}{\sim} -\frac{1}{64}\frac{z^3(z^3-8)^3}{(z^3+1)^3},\leqno{17.}$$  
$$\nu(\f O_1)=\{2,2,2,2\}, \ \ \ \ \nu(\f O_2)=\{2,3,6\},$$  
\vskip 0.1cm
\noindent In particular, if $\deg A>12$, then $\f O_2=\f O_2^A$, $\f O_1=\f O_1^A$.
\et
\pr It follows from  \eqref{elki+} that $\chi(\f O_2^A)\geq \chi(\f O_2)$ and the equality is attained if and only if $\f O_2^A= \f O_2$. Therefore,  if $\chi(\f O_2^A)= 0$, then   $\f O_2^A=\f O_2$ and hence 
$\f O_1^A=\f O_1$ since \eqref{us} implies that for any covering map $A:\, \f O_1\rightarrow \f O_2$ the orbifold $\f O_1$ is defined by $\f O_2$ in a unique way.
So, below we will assume that  $\chi(\f O_2^A)> 0.$ 

We will denote by $\nu_A$ the ramification function of $\f O_2^A$ and by $\nu$   the ramification function of
 $\f O_2$. We also  
will use the notation 
$$R(f)=\Big(\{l_{11},l_{12}, \dots l_{1s_1}\}_{z_1},  \dots,  \{l_{r1},l_{r2}, \dots, l_{rs_r}\}_{z_r}\Big)$$
for denoting that a rational function $f$  has $r$ critical values $z_1,$ $z_2,$ ... $z_r$, and 
the collection of local degrees of $f$ at points of the set $f^{-1}\{z_i\}$, $1\leq i \leq r,$  is $\{l_{i1},l_{i2}, \dots, l_{is_i}\}.$

As in Lemma \ref{lemlas}, the conditions  $\f O_2^A\preceq \f O_2$ and $\chi(\f O_2)=0$, $\chi(\f O_2^A)>0$ impose strong restrictions on possible collections $\nu(\f O_2^A)$, and an easy ana\-lysis of lists \eqref{list} and  \eqref{list2} shows that   either 
$\nu(\f O_2^A)=\{n,n\},$ $n\leq 4,$  or 
 $\nu(\f O_2^A)=\{2,2,n\},$ $n\leq 4,$ or $\nu(\f O_2^A)=\{2,3,3\}.$

\vskip 0.2cm

\noindent {\bf Case 1:} $\boldsymbol{\nu(\f O_2^A)=\{n,n\}.}$ 
If $n=2$, then  $\f O_2^A\preceq \f O_2$ implies that $\nu(\f O_2)$ is either $\{2,2,2,2\}$, or $\{2,4,4\}$, or  $\{2,3,6\}$.
Assume say that $\nu(\f O_2)=\{2,4,4\}$ and let 
$x_1,x_2,y_1,y_2,y_3\in \C\P^1$ be points such that 
\be \l{pedf} \nu_A(x_1)=2, \ \ \ \nu_A(x_2)=2,\ee
\be \l{zr} \nu(y_1)=2, \ \ \ \nu(y_2)=4,\ \ \ \nu(y_3)=4.\ee
Then either
$$\{x_1,x_2\}=\{y_1,y_2\},$$ or
$$\{x_1,x_2\}=\{y_1,y_3\},$$ or
\be \l{vr1} \{x_1,x_2\}=\{y_2,y_3\}.\ee
Further, since  \eqref{pedf} implies that 
$$\f R(A)=\left(\{2\}_{x_1}, \{2\}_{x_2}\right),$$ it follows from  \eqref{us} 
that in the first two cases 
$$\nu(\f O_1)=\{2,4,4\},$$
while in the third one  
$$\nu(\f O_1)=\{2,2,2,2\}.$$ 
Thus, we arrive to cases 2 and 3 listed in the theorem.

Similarly, if  $\nu(\f O_2)=\{2,3,6\}$ and 
$y_1,y_2,y_3\in \C\P^1$ are points such that 
$$\nu(y_1)=2, \ \ \ \nu(y_2)=3,\ \ \ \nu(y_3)=6,$$
then  $$\{x_1,x_2\}=\{y_1,y_3\}$$ and  
$$\nu(\f O_1)=\{3,3,3\}.$$ 
Finally, if $\nu(\f O_2)=\{2,2,2,2\}$ we conclude that  
$$\nu(\f O_1)=\{2,2,2,2\}.$$

The cases $n=3$ and $n=4$ are considered in the same way as above. Namely, 
if  $\nu(\f O_2^A)=\{3,3\},$ then $\nu(\f O_2)$ is either $\{3,3,3\}$, or  $\{2,3,6\}$, and we arrive to cases 5 and 6, correspondingly, while if $\nu(\f O_2^A)=\{4,4\},$ then $\nu(\f O_2)=\{2,4,4\}$ and we arrive to case 7.

\vskip 0.2cm

\noindent {\bf Case 2:} $\boldsymbol{\nu(\f O_2^A)=\{2,2,n\}.}$ The proof goes as above with some modifications. 
Let 
$x_1,x_2,x_3\in \C\P^1$ be points such that 
$$\nu_A(x_1)=2, \ \ \ \nu_A(x_2)=2, \ \ \ \nu_A(x_3)=n.$$
Assume say that $n=3$. Then 
$\nu(\f O_2)=\{2,3,6\}$, and if  
$y_1,y_2,y_3\in \C\P^1$ are points such that 
$$\nu(y_1)=2, \ \ \ \nu(y_2)=3,\ \ \ \nu(y_3)=6,$$ then 
$$\{x_1,x_2\}=\{y_1,y_3\},\ \  \ x_3=y_2.$$ 
Now however we must consider two types of branching of $A$ corresponding to the possibilities 
$A\s \frac{1}{2}(z^3+z^{-3})$ and $A\s T_3.$ In the first case
$$\f R(A)=\Big (\{2,2,2\}_{x_1}, \{2,2,2\}_{x_2}, \{3,3\}_{x_3}\Big),$$
in the second
$$\f R(A)=\Big (\{1,2\}_{x_1}, \{1,2\}_{x_2}, \{3\}_{x_3}\Big).$$
Correspondingly, either 
$$\nu(\f O_1)=\{3,3,3\}.$$
or 
$$\nu(\f O_1)=\{2,3,6\}.$$

Similarly, if $n=4$, then  
$\nu(\f O_2)=\{2,4,4\}$, and either 
$$\f R(A)=\Big (\{2,2,2,2\}_{x_1}, \{2,2,2,2\}_{x_2}, \{4,4\}_{x_3}\Big),$$
or
\be \l{xrun} \f R(A)=\Big (\{1,1,2\}_{x_1}, \{2,2\}_{x_2}, \{4\}_{x_3}\Big),\ee
or 
\be \l{xrun2} \f R(A)=\Big (\{2,2\}_{x_1}, \{1,1,2\}_{x_2}, \{4\}_{x_3}\Big).\ee
In the first case, $$\nu(\f O_1)=\{2,2,2,2\},$$ while in each of cases \eqref{xrun} and \eqref{xrun2}, 
either \be \l{ili1} \nu(\f O_1)=\{2,2,2,2\},\ee or \be \l{ili2} \nu(\f O_1)=\{2,4,4\}.\ee Say, if \eqref{xrun} holds, and 
$y_1,y_2,y_3\in \C\P^1$ are the points such that 
$$\nu(y_1)=2, \ \ \ \nu(y_2)=4,\ \ \ \nu(y_3)=4,$$  then 
\eqref{ili1} holds if $x_1=y_1,$ while \eqref{ili2} holds if $x_1=y_2$ or $x_1=y_3$. 

Finally, if $n=2$, then $A\s \frac{1}{2}(z^2+z^{-2}),$ and  $\nu(\f O_2)$ is 
either 
$\{2,2,2,2\}$ or $\{2,4,4\}$.
In the both  cases,  $$\nu(\f O_1)=\{2,2,2,2\}.$$

\vskip 0.2cm

\noindent {\bf Case 3:} $\boldsymbol{\nu(\f O_2^A)=\{2,3,3\}.}$ 
In this case $\nu(\f O_2)=\{2,3,6\}$, and considering three possible branching type for tetrahedral functions 
\be \l{pizd2}\f R(A)=\Big (\{2,2\}_{z_1}, \{1,3\}_{z_2}, \{1,3\}_{z_3}\Big),\ee
\be \l{pizd3} \f R(A)=\Big(\{1,1,2,2\}_{z_1}, \{3,3\}_{z_2}, \{3,3\}_{z_3}\Big), \ee
\be \l{pizd1} \f R(A)=\Big(\{2,2,2,2,2,2\}_{z_1}, \{3,3,3,3\}_{z_2}, \{3,3,3,3\}_{z_3}\Big),
\ee
we arrive to the cases  15, 16, 17 correspondingly. \qed

\br Modifying the above proof  one can see that all the possibilities listed in Theorem \ref{l} actually occur. 
For example,  
for any rational function $A\sim z^2$,  {\it there exist} orbifolds $\f O_1$ and $\f O_2$ such that $\nu(\f O_1)=\{2,2,2,2\}$, $\nu(\f O_2)=\{2,4,4\}$,
and  $A:\, \f O_1\rightarrow \f O_2$ is a covering map. Indeed, let $x_1,$ $x_2$ be points such that \eqref{pedf} holds. 
Define $\f O_2$ by formula  \eqref{zr}, where $y_2,y_3$ satisfy  \eqref{vr1} and $y_1$ is taken arbitrary,
and then define $\f O_1$ by formula \eqref{us}. 
Since $\f O_2^A\preceq \f O_2$ implies that for any $z\in \C\P^1$ the number $\deg_z A$ divides the number $\nu(A(z))$, the orbifold $\f O_1$ is well-defined and  $A:\, \f O_1\rightarrow \f O_2$ is a covering map. 
Thus,  Theorem \ref{l} gives  
a complete list of $\mu$-equivalence classes of rational functions $A$ which fit diagram \eqref{xxuu} but satisfy $g(\t S_A)=0$ instead of $g(\t S_A)=1.$
\er

\bc \l{ll+} Let $A$ be a rational function and $\f O_1$, $\f O_2$ orbifolds such that $\nu(\f O_1)=\nu(\f O_2)$ and 
$A:\, \f O_1\rightarrow \f O_2$ is a covering map between orbifolds. 
Then 
either  $\chi(\f O_2^A)=0$ and $\f O_2=\f O_2^A$, $\f O_1=\f O_1^A$, or $\chi(\f O_2^A)>0$ and one of the following conditions holds:

$$\nu(\f O_2^A)=\{2,2\},  \ \   {A}\underset{\mu}{\sim}z^2, \ \ \nu(\f O_1)=\nu(\f O_2)=\{2,2,2,2\},\leqno{1.}$$ 
$$\nu(\f O_2^A)=\{2,2\},  \ \   {A}\underset{\mu}{\sim}z^2, \ \ \nu(\f O_1)=\nu(\f O_2)=\{2,4,4\},\leqno{2.}$$ 
$$\nu(\f O_2^A)=\{3,3\},  \ \   {A}\underset{\mu}{\sim}z^3, \ \ \nu(\f O_1)=\nu(\f O_2)=\{3,3,3\},\leqno{3.}$$  
$$\nu(\f O_2^A)=\{2,2,2\}, \ \   {A}\underset{\mu}{\sim} \frac{1}{2}\left(z^2+\frac{1}{z^2}\right),\ \  \nu(\f O_1)=\nu(\f O_2)=\{2,2,2,2\},\leqno{4.}$$ 
$$\nu(\f O_2^A)=\{2,2,3\}, \ \    {A}\underset{\mu}{\sim} T_3,  \ \  \nu(\f O_1)=\nu(\f O_2)=\{2,3,6\},\leqno{5.}$$  
$$\nu(\f O_2^A)=\{2,2,4\}, \ \ \   {A}\underset{\mu}{\sim} T_4,  \ \ \ \nu(\f O_1)=\nu(\f O_2)=\{2,4,4\},\leqno{6.}$$ 
$$\nu(\f O_2^A)=\{2,3,3\}, \ \ {A}\underset{\mu}{\sim} -\frac{1}{64}\frac{z(z-8)^3}{(z+1)^3},\ \ \nu(\f O_1)=\nu(\f O_2)=\{2,3,6\}.\leqno{7.} $$ 
In particular, if $\deg A>4$, then $\f O_2=\f O_2^A$, $\f O_1=\f O_1^A$.
\ec
\pr The corollary follows from Theorem \ref{l} since $\nu(\f O_1)=\nu(\f O_2)$ implies the equality $\chi(\f O_1)=\chi(\f O_2)=0$ by \eqref{rhor+}. \qed

\vskip 0.2cm

Recall that Latt\`es maps are rational functions which can be defined in one of the following ways (see \cite{mil2}, \cite{lattes}). 
First, a Latt\`es map $A$ may defined by the condition that
there exist a Riemann surface 
$\f C$ of genus one and holomorphic maps $\alpha: \f C\rightarrow \f C$ and $\pi:\f C\rightarrow \C\P^1$
such that the diagram 
\be \l{xxuuii}
\begin{CD}
\f C @>\alpha>> \f C \\
@VV\pi V @VV\pi V\\ 
\C\P^1 @>A >> \ \ \C \P^1
\end{CD}
\ee
is commutative. This condition is equivalent to the apparently stronger condition that $\pi$ in \eqref{xxuuii} is a quotient map
$\pi: {\f C}\rightarrow   
{\f C}/\Omega$ for some finite subgroup  $\Omega\subseteq {\rm Aut}({\f  C}).$ 
Finally, a Latt\`es map $A$ may defined by the condition that there exists an orbifold $\f O$ in $\C\P^1$ such that $\chi(\f O)=0$ and 
$A:\f O\rightarrow \f O$ is a covering map between orbifolds. Thus, Corollary \ref{ll+} implies the following corollary.

\bc \l{ll} For any Latt\`es map $A$ of degree greater than four the equality   $g(\t S_A)=1$ holds. \qed
\ec

\br It is easy to see that there exist rational functions $A$ with $g(\t S_A)=1$ which are not Latt\`es maps. 
Indeed, let  $A:\f O_1\rightarrow \f O_2$ be any covering map  between orbifolds such that $\f O_1\neq \f O_2$ and $\deg A>12.$
Then it follows from Theorem  \ref{mt3} that  
$\chi(\f O_2^A)=0$ and $\f O_2=\f O_2^A$, $\f O_1=\f O_1^A$. Thus,   $g(\t S_A)=1$ by Lemma \ref{ml}. On the other hand,  if 
$\f O$ is an orbifold  such that 
$A:\f O\rightarrow \f O$ is a covering map between orbifolds, then 
it follows from  \eqref{elki+} that $$\chi(\f O_2^A)\geq \chi(\f O),\ \ \ \ \ \  \chi(\f O_1^A)\geq \chi(\f O)$$ and the equality is attained if and only if $\f O_2^A= \f O$, $\f O_1^A= \f O.$ 
Since 
$$\chi(\f O_2^A)=\chi(\f O_1^A)=\chi(\f O)=0,$$ this implies that 
$\f O_2=\f O_1=\f O,$ in contradiction with  $\f O_1\neq \f O_2$.
\er


\end{section}

\bibliographystyle{amsplain}

\providecommand{\bysame}{\leavevmode\hbox to3em{\hrulefill}\thinspace}
\providecommand{\MR}{\relax\ifhmode\unskip\space\fi MR }
\providecommand{\MRhref}[2]{%
  \href{http://www.ams.org/mathscinet-getitem?mr=#1}{#2}
}
\providecommand{\href}[2]{#2}
\begin{thebibliography}{}

\end{thebibliography}


\begin{thebibliography}{10}

\bibitem {az} R. Avanzi, U. Zannier, \textit{
The equation $f(X)=f(Y)$ in rational functions $X=X(t),$ $Y=Y(t)$,}
Compositio Math. 139 (2003), no. 3, 263-295. 

\bibitem {bilu} Y. Bilu, R. Tichy, 
\textit{The Diophantine equation $f(x) = g(y)$},
Acta Arith. 95, No.3, 261-288 (2000).


\bibitem{couv} J.-M. Couveignes, \textit{Calcul et rationalit\'e de fonctions de Belyi en genre 0}, Ann. de l'Inst. Fourier 44 (1994) 1-38.


\bibitem{dan}  L. Danilov, \textit{The Diophantine equation $x^3-y^2=k$ and a conjecture of M. Hall},  Mat. Zametki 32 (1982), no. 3, 273-275.



\bibitem{e2} A. Eremenko, \textit{Some functional equations connected with the iteration of rational functions} (Russian), Algebra i Analiz 1 (1989), 102-116; translation in Leningrad Math. J. 1 (1990), 905-919. 


\bibitem{e} A. Eremenko, \textit{Invariant curves and semiconjugacies of rational functions},  Fundamenta Math., 
219, 3 (2012) 263-270.

\bibitem {fk} H. Farkas, I. Kra,  \textit{Riemann surfaces}, Graduate Texts in Mathematics, 71. Springer-Verlag, New York, 1992.

\bibitem {fsch} M. Fried, 
\textit{
On a conjecture of Schur},
Michigan Math. J. 17, 1970, 41-55. 

\bibitem {f1} M. Fried, 
\textit{On a theorem of Ritt and related diophantine problems}, 
J. Reine Angew. Math. 264, 40-55 (1973).




\bibitem {fried} M. Fried, {\it Introduction to modular towers: generalizing dihedral group-modular curve connections,} Recent developments in the inverse Galois problem, 111-171, Contemp. Math., 186, Amer. Math. Soc., Providence, RI, 1995.
 

\bibitem {des1}  E. Girondo,  G. Gonz\'alez-Diez,  \textit{Introduction to compact Riemann surfaces and dessins d'enfants,} London Mathematical Society Student Texts, 79. Cambridge University Press, Cambridge, 2012.



\bibitem {klein} F. Klein, 
\textit{Lectures on the icosahedron and the solution of equations of the fifth degree},
New York: Dover Publications, (1956).





\bibitem {kt} D. Kreso, R. Tichy, 
\textit{Diophantine equations in separated variables}, Period Math Hung., to appear.


\bibitem {lazv} S. Lando and A.K. Zvonkin,  \textit{Graphs on surfaces and their applications,} Springer-Verlag, 2004.


\bibitem {mazv} N. Magot, A. Zvonkin,  \textit{Belyi functions for archimedean solids, Discrete Math.} 217 (2000) 249-271.



\bibitem {ms}  A. Medvedev, T. Scanlon,  \textit{
Invariant varieties for polynomial dynamical systems,}  Annals of Mathematics, 179 (2014), no. 1, 81 - 177. 



\bibitem{mil} J. Milnor, \textit{Dynamics in one complex variable}, Princeton Annals in Mathematics 160. Princeton, NJ: Princeton University Press (2006).

\bibitem{mil2} J. Milnor, \textit{On Latt\`es maps,} Dynamics on the Riemann Sphere. Eds. P. Hjorth and C. L. Petersen. A Bodil Branner Festschrift, European Mathematical Society, 2006, pp. 9-43.






\bibitem {mp2} M. Muzychuk, F. Pakovich, {\it  Jordan-Holder theorem for imprimitivity systems and maximal decompositions of rational functions}, 
Proc. Lond. Math. Soc., 102 (2011) , no. 1, 1-24.



 \bibitem {prime} F. Pakovich, {\it  Prime and composite Laurent polynomials, Bull. Sci. Math.,} 133 (2009), 693-732.



\bibitem {semi} F. Pakovich, {\it On semiconjugate rational functions,} Geom. Funct. Anal., 26 (2016), 1217-1243. 
 
 

\bibitem {pj} F. Pakovich, {\it    Polynomial semiconjugacies, decompositions of iterations, and invariant curves},  Ann. Sc. Norm. Super. Pisa Cl. Sci., to appear. 




 \bibitem {alc} F. Pakovich, {\it       Algebraic curves $P(x)-Q(y)=0$ and functional equations,}  Complex Var. and Elliptic Equ., 56 (2011), no. 1-4, 199-213.


\bibitem {amer} F. Pakovich, {\it    On the equation $P(f)=Q(g),$ where $P,Q$ are polynomials and $f,g$ are entire functions,} Amer. Journal of Math., 132 (2010), no. 6, 1591-1607.


\bibitem {gen} F. Pakovich, {\it  On algebraic curves $A(x)-B(y)=0$ of genus zero}, Math. Z., to appear. 



\bibitem {pz} F. Pakovich, A. Zvonkine, {\it   Davenport-Zannier  Polynomials over $\Q$},  	Int. J. Number Theory, to appear.


\bibitem {lattes} F. Pakovich, {\it On generalized Latt\`es maps}, preprint, arxiv:1612.01315.




\bibitem{r2} J. F. Ritt, \textit{Prime and composite polynomials,} Trans. Amer. Math. Soc., 23, 1922, 51-66.



\bibitem{wi} H. Wielandt,
\textit{Finite Permutation Groups},
 Academic Press, 1964, Berlin.



\end{thebibliography}

\end{document}